\documentclass[11pt]{amsart}
\newsymbol\subsetneqq 2324
\newtheorem{thm}[equation]{Theorem}
\newtheorem{pro}[equation]{Proposition}
\newtheorem{cor}[equation]{Corollary}
\newtheorem{lem}[equation]{Lemma}
\newtheorem{exa}[equation]{Example}

\newtheorem{DEF}[equation]{Definition}
\newtheorem{rem}[equation]{Remark}

\numberwithin{equation}{section}

\def\a1s{a_1,\cdots, a_s}
\def\a{\alpha}

\def\aa{\mathcal A}

\def\adot{\dot{\alpha}}

\def\andd{\quad\hbox{and}\quad}

\def\ad{\hbox{ad}}

\def\b{\beta}

\def\bb{\mathcal{B}}

\def\bl4{B_{\ell\geq4}}

\def\bbbc{{\mathbb C}}

\def\d{\delta}
\def\D{\Delta}

\def\dd{\mathcal D}

\def\ve{\varepsilon}

\def\bbbf{\mathbb{F}}

\def\GL{GL}
\def\gg{{\mathcal G}}
\def\dg{\dot{\gg}}
\def\dgg{\dot{\gg}}

\def\fg{\mathfrak{g}}

\def\hh{{\mathcal H}}
\def\dhh{\dot{\hh}}
\def\hhz{\hh^0}
\def\sdhh{\dhh^\star}

\def\hhstar{\hh^{\star}}

\def\fh{\mathfrak{h}}

\def\bbbk{\mathbb{K}}
\def\lam{\lambda}
\def\Lam{\Lambda}
\def\LL{\mathcal{L}}
\def\LLc{\mathcal{L}_c}
\def\LLa{\mathcal{L}_{\a}}
\def\mod{\hbox{mod}}
\def\bigop{\bigoplus}
\def\ep{\epsilon}
\def\fm{(\cdot,\cdot)}
\def\bbbq{\mathbb{Q}}
\def\bbbn{\mathbb{N}}

\def\bbbr{{\mathbb R}}

\def\rd{\dot{R}}

\def\rds{\dot{R}_{sh}}
\def\rdl{\dot{R}_{lg}}
\def\rde{\dot{R}_{ex}}

\def\tr{\hbox{tr}}

\def\1k{\frac{1}{k}}
\def\op{\oplus}
\def\ot{\otimes}

\def\la{\langle}
\def\ra{\rangle}

\def\qed{\hfill$\Box$}
\def\sub{\subseteq}
\def\sg{\sigma}

\def\rtimes{R^{\times}}

\def\rzero{R^{0}}

\def\hd{\dot{\hh}}

\def\tr{\hbox{tr}}

\def\T{{\mathcal T}}

\def\u{{\mathcal U}}
\def\ud{\dot{\mathcal U}}

\def\v{{\mathcal V}}

\def\vz{{\mathcal V}^{0}}

\def\zn{{\mathbb Z}^{\nu}}

\def\bbbz{{\mathbb Z}}

\def\1il{1\leq i\leq\ell}

\begin{document}
\markboth{A generalization of extended affine Lie algebras}{}
\title{ A generalization of extended affine Lie algebras}
\maketitle
 \setcounter{section}{-1}
 \vspace{1cm} \noindent
\vspace{.5cm} \begin{center}Malihe Yousofzadeh\footnote{
Department of Mathematics, University of Isfahan, Isfahan, Iran,
P. O. Box 81746-73441, E-mail: ma.yousofzadeh@math.ui.ac.ir}
\end{center}

ABSTRACT. We introduce a new class of possibly infinite
dimensional Lie algebras and study  their structural properties.
Examples of this new class of Lie algebras are finite dimensional
 simple Lie algebras containing a nonzero split torus, affine and extended affine Lie algebras.
 Our results generalize well-known properties  of these examples.
\section{ Introduction}
In 1985,  Saito  [Sa] introduced the notion of  an {\it extended
affine root system} which is a nonempty subset $R$ of  a finite
dimensional real vector space $V,$ equipped with  a positive
semidefinite bilinear form $(\cdot,\cdot)$, satisfying certain
axioms. Following [Az2], we call them  {\it Saito's extended
affine root systems} or SEARSs for short. In 1990,
H$\phi$egh-Krohn and  Torresani [H-KT] introduced a new class of
Lie algebras over the field of complex numbers. The basic features
of these   Lie algebras are the existence of a non-degenerate
symmetric invariant bilinear form, a finite dimensional Cartan
subalgebra, a discrete root system and the ad-nilpotency of the
root spaces attached to non-isotropic roots.
In 1997,  Allison,   Azam,   Berman,   Gao and  Pianzola [AABGP]
called these Lie algebras  {\it extended affine Lie algebras} (or
EALAs for short). The subalgebra of an EALA generated by its
non-isotropic root spaces, called the {\it core}, plays a very
important role in the study of EALAs. Namely, up to the choice of
a certain derivation and a 2-cocycle,  an EALA is determined by
its core modulo center, called the {\it centerless core} [N1,2].
In [BGKN], [BGK], [AG] and [Yo1], the authors give a complete
description of the centerless cores of reduced extended affine Lie
algebras (see also [Az1]).

In [AABGP], the authors take the properties of the root system
 of an EALA as  axioms for  the new notion of an {\it extended affine root systems} or EARS for short.
The non-isotropic roots  of an  EARS form  a SEARS, but not all
SEARS arise in this way. In fact,  there is a one to one
correspondence between reduced SEARS and non-isotropic parts of
EARSs [Az2].

There had been  several attempts  to construct  classes of  Lie
algebras  with SEARSs as their root systems [Sa, Ya1, SaY, Ya2].
But the sets of non-isotropic roots corresponding to  all such
constructed Lie algebras are reduced.

In 2004,  Yoshii [Yo2] introduced the notion of a {\it
$(\D,G)-$graded }Lie algebra, $\Delta$ a finite irreducible root
system and $G$ an abelian group.  These provide examples of Lie
algebras with not necessarily reduced SEARSs as their root
systems. Also in [AKY], the authors axiomatically introduced a
class   of possibly infinite dimensional Lie algebras over a field
of characteristic zero called {\it toral type extended affine Lie
algebras}. The non-isotropic
 roots  of a toral type extended affine Lie algebra form a not necessarily reduced  SEARS.  Roughly speaking, they replace the existence of a Cartan subalgebra in the axioms of an
EALA with the existence of a nonzero split torus. This then allows
the existence of divisable roots for the corresponding root
systems.

We  introduce a class $\T$ of  Lie algebras whose axioms are  a
modified version of axioms introduced in [AKY].  The main purpose
of this work is to study the Lie algebras in the class $\T$ along
the  lines one has studied EALAs. In Section $1$, we introduce the
axioms for the Lie algebras in $\T$ and investigate the basic
properties of the algebras.
 In Section 2, we study
 the core $\LLc$ of a Lie algebra $\LL\in\T$. We prove that
 $\LLc$ is a $(\D, G)-$graded Lie algebra. We also describe a {\it loop}
construction of Lie  algebras in $\T$ which satisfy a
 {\it division} property (see subsection 2.2). We conclude the section with an example  of
a  Lie  algebra in $\T$ of type $C$ and of
  arbitrary nullity, with non-isotropic root spaces
having  dimension greater than $1$, a phenomenon which does not
happen for EALAs. In Section 3, we study the elements in $\T$ of
nullity zero. We find a necessary and sufficient condition for   a
Lie algebra  to be an element  of $ \T$  of nullity zero. We also
show that the  class $\T$ is a generalization of the class of
finite dimensional simple Lie algebras containing a nonzero split
torus. In the appendix we  recall the definition of a toral type
extended affine Lie algebra. We prove that  a Lie algebra in $\T$
satisfies the axioms of a  toral type extended affine Lie algebra.
\section{  A new class of infinite dimensional Lie algebras}
Throughout   this work $\bbbf$ is a field of characteristic zero.
Unless otherwise  mentioned, all vector spaces are considered over
$\bbbf$. In the present  paper,  we denote the dual space  of a
vector space $V$ by $V^\star$.  If a finite dimensional vector
space $V$ is equipped with a non-degenerate symmetric bilinear
form, then for $\a\in V^\star$, we take $t_\a$ to be the unique
element in $V$ representing $\a$ through the form. Also for an
algebra $\aa$, $Z(\aa)$ denotes the center of $\aa.$ We refer to a
finite root system as a subset  $\D$ of a vector space so that
$0\in\D$ and $\D\setminus\{0\}$ is a finite root system in the
sense of [Bo]. For a finite root system $\D$, we set
$\D^\times:=\D\setminus\{0\}$ and
$\D_{red}:=\{0\}\cup\{\a\in\D\mid\a/2\not\in\D\}$. Also we make
the convention that  for a finite root system all roots are short
if only one root length occurs.

\begin{DEF}
{\em Let $\LL$ be a Lie algebra.  A nontrivial subalgebra $\hh$ of
$\LL$ is called a {\it  split toral subalgebra} or  a {\it  split
torus} if
$$\LL=\bigoplus_{\a\in\hh^\star}\LL_\a,\hbox{ where
}\LL_\a=\{x\in\LL\mid[h,x]=\a(h)x\hbox{ for all }h\in\hh\}.$$One
can see that a split toral subalgebra  is abelian.
 }
\end{DEF}
As all toral subalgebras [Se, \S I.1] occurring in this paper are
split, in the sequel, in place of saying  split toral subalgebra,
we say toral subalgebra.

Now let $\hh$ be a finite dimensional   toral subalgebra of a Lie
algebra  $\LL=\bigoplus_{\a\in\hh^\star}\LL_\a$  equipped with a
non-degenerate symmetric invariant bilinear form $(\cdot,\cdot)$
such that the form $(\cdot,\cdot)$ restricted to  $\hh$ is
non-degenerate. Then the Jacobi identity implies that
$[\LL_\a,\LL_\b]\subseteq\LL_{\a+\b}$ for all $\a,\b\in\hh^\star$
and the invariance of the form implies that for
$\a,\b\in\hh^\star$,
\begin{equation} \label{v4}
(\LL_\a,\LL_\b)=0 \;\;\;\hbox{unless}\;\;\a+\b=0.
\end{equation}
We call $R:=\{\a\in\hh^\star\mid\LL_\a\not=\{0\}\}$, the root
system  of $\LL$ with respect to $\hh$.  We will transfer the form
from $\hh$ to $\hh^\star$ by requiring $(\a,\b):=(t_\a,t_\b)$ for
any $\a,\b\in\hh^\star$. This allows us to define $$
\rtimes:=\{\a\in R\mid (\a,\a)\not=0\}\andd R^0:=\{\a\in R\mid
(\a,\a)=0\}.$$
 The elements of $\rtimes$ and of $R^0$ are called {\it
 non-isotropic}
and {\it isotropic} roots of $R$.
\begin{DEF}\label{core} {\em
The subalgebra $\LLc$ of $\LL$ generated by the root spaces
$\LL_\a$, $\a\in\rtimes$, is called the {\it core} of $\LL$.
Following [AY, \S 5], we call $\LL_{cc}:={\LLc}/{Z(\LLc)}$ the
{\it centerless core } of $\LL$.  The algebra  $\LL$ is called
{\it tame} if $\LL_c$ contains its centralizer in $\LL$, namely
$C_\LL(\LL_c)\sub \LL_c$. }
\end{DEF}
It follows from the finite dimensionality of $\hh$ and the
non-degeneracy of the form on $\hh$ and on $\LL_0$ that
$\LL_0=\hh\op\hh^\bot,$ where $\hh^\bot=\{x\in \LL_0\mid (h,x)=0
\;\hbox{for all }\; h\in\hh\}.$ So for $\a\in R$, $x\in\LLa$ and
$y\in\LL_{-\a},$ we have $[x,y]=t_{x,y}+h_{x,y}$ for some
$t_{x,y}\in\hh\hbox{\;and\; } h_{x,y}\in \hh^\bot$.  Put
\begin{equation}
\hh_\a:=\hbox{span}_\bbbf\{h_{x,y}\mid
x\in\LLa,y\in\LL_{-\a}\},\;\a\in R.\label{new}
\end{equation} Now let $\a\in
R$, $x\in\LLa$, $y\in\LL_{-\a}$ and $h\in \hh$. Then
\begin{eqnarray*}
(h,t_{x,y})=(h,[x,y])=([h,x],y)=\a(h)(x,y)&=&(t_\a,h)(x,y)\\
&=&((x,y)t_\a,h).
\end{eqnarray*}
Since the form  is symmetric and  non-degenerate on $\hh$,  we
conclude that  $t_{x,y}=(x,y)t_\a$. So
\begin{equation}\label{m0}
\begin{array}{c}
[x,y]=(x,y)t_\a+h_{x,y}; \;\a\in R,\; x\in \LL_\a,\;
y\in\LL_{-\a}\;\;\hbox{\;and\;}
\\
\hbox{if \;} 0\neq\a,\;x\in\LL_\a,\;y\in\LL_{-\a}\andd [x,y]=t_\a,
\hbox{\;then\;}(x,y)=1.
\end{array}
\end{equation}

\begin{DEF}
{\em We denote by $\T$ the class   of all triples  $(\LL,
(\cdot,\cdot),\hh)$ satisfying the following six axioms:

{\bf (T1)} $\LL$ is a Lie algebra over $\bbbf$ equipped with a
non-degenerate symmetric invariant bilinear form $\fm$.

{\bf (T2)} $\hh$ is a finite dimensional  toral subalgebra of
$\LL$ such that $\fm|_{\hh\times \hh}$ is non-degenerate. Let $R$
be the root system of $\LL$  with respect to $\hh.$

{\bf (T3)} For $\d\in R^0$, there exist $x\in\LL_\d$ and
$y\in\LL_{-\d}$ such that $[x,y]=t_\d$. Also  for each $\a\in
\rtimes$ and $0\neq x\in\LL_\a$, there exists $y\in \LL_{-\a}$
such that $[x,y]=t_\a.$

We point out  that the conditions on $\rtimes$ and on $R^0$ are
not symmetric. More precisely (T3) requires that every $0\neq
x\in\LL_\a,$ $\a\in\rtimes,$ is part of an $\frak{sl_2}-$triple
while  for $\a\in R^0$ only the existence of certain elements is
required.

 {\bf (T4)} If $\a\in\rtimes$ and
$x_\a\in\LL_\a$, then $\ad_{\LL}x_\a$ acts locally nilpotently on
$\LL$.

{\bf (T5)} $R$ is irreducible, in the sense that it satisfies the
following two conditions:

(a) $\rtimes$ cannot  be written as a disjoint union of two
orthogonal nonempty subsets of $\rtimes.$

(b) For $\d\in R^0$ there exists $\a\in\rtimes$ such that
$\a+\d\in R$.

{\bf (T6)} $\Lam:=\la R^0\ra=\hbox{span}_{\bbbz}(R^0)$ is a free
abelian group of finite rank (we consider $\{0\}$ as a free
abelian group of rank  $0$).

When $(\cdot,\cdot)$ and $\hh$ are fixed, we denote a triple
$(\LL, (\cdot,\cdot),\hh)\in \T$ by $\LL.$ For $\LL\in\T$ with
root system $R$, the rank of $\Lam=\la R^0\ra$ is called the {\it
nullity} of $\LL$. }
\end{DEF}
Axioms $(T1)-(T6)$ are a modified version of axioms of the class
of {\it toral type extended affine Lie algebras} introduced in
[AKY].
\\
$ $

 Now let $(\LL,(\cdot,\cdot),\hh)$ satisfy $(T1)-(T4)$. Then  (\ref{m0})
together with axiom ($T3$) gives:
\begin{equation}\label{form1}
[\LLa,\LL_{-\a}]=\bbbf t_{\a}\op\hh_\a,\;\;\a\in R
\end{equation}
For  $\a\in\rtimes,$ set $h_\a:=2t_\a/(t_\a,t_\a)$. Then  by
$(T3)$, there exist $e_{\pm\a}\in\LL_{\pm\a}$ such that
$(e_\a,h_\a,e_{-\a})$ is an $\mathfrak{sl}_2\hbox{-triple},$ that
is $[h_\a,e_\a]=2e_\a,[h_\a,e_{-\a}]=-2e_{-\a}$ and
$h_\a=[e_\a,e_{-\a}]$. Now let  $\a\in\rtimes$. The reflections
$w_\a\in\GL(\hh^\star)$ and $\check{w}_\a\in\GL(\hh)$ are defined
by
$$\begin{array}{cl}
 w_\a(\b)=\b-\frac{2(\b,\a)}{(\a,\a)}\a,&\b\in\hh^\star\\
 \check{w}_\a(h)=h-\frac{2(h,t_\a)}{(t_\a,t_\a)}t_\a,& h\in\hh.
 \end{array} $$
Using $(T4)$, for $t\in \bbbf\setminus\{0\}$ we can define an
automorphism of $\LL$ by $\hbox{exp}(\ad
te_\a):=\sum_{k=0}^\infty{((\ad te_\a)^k}/{k!}). $ So for  $
\a\in\rtimes$ and $t\in\mathbb{F}\setminus\{0\},$ we have
$$\theta_{\a}(t):=\hbox{exp}(\ad(te_{\a}))\hbox{exp}(\ad(-t^{-1}e_{-\a}))\hbox{exp}(\ad(te_{\a}))\in\hbox{Aut}(\LL).$$
One  can easily check that for $\a\in\rtimes$, $h\in\hh$ and
$t\in\bbbf\setminus\{0\},$
\begin{equation} \theta_\a(t)(h)=\check{w}_\a(h). \label{v1}
\end{equation}
\begin{pro}\label{v2}
Let $\LL$ satisfy $(T1)-(T4)$. Then for  $\a\in\rtimes$ and $\b\in
R$, we have

$(i)$ $ \theta_{\a}(t)(\LL_{\b})=\LL_{w_\a(\b)}$.

$(ii)$ $ 2(\b,\a)/(\a,\a)\in\bbbz.$

$(iii)$ $\dim(\LL_\a)=\dim(\LL_{-\a})$.
\end{pro}
\noindent{\bf Proof.}  Parts $(i)$ and $(ii)$ can be seen using
the same arguments as in  [MP, Proposition 4.1.4] and  [AABGP,
Proposition I.1.29] respectively.

 $(iii)$ By  part ($i$),   we have
$\theta_{\a}(1)(\LL_{\a})=\LL_{w_\a(\a)}=\LL_{-\a}$. So
$\dim(\LL_\a)=\dim(\LL_{-\a}).$\qed
\\
$ $

For the further study of the  algebras in $\T$, we need to assume
one more axiom. Let $(\LL,(\cdot,\cdot),\hh)$ satisfy $(T1)-(T5)$.
Since $0\in R^0$, $(T5)(b)$ implies that  $\rtimes \neq\emptyset.$
So there exist $\a\in\rtimes $  and $k\in\bbbf$ such that
$k(\a,\a)=1$. By rescaling the form if necessary,  we may assume
that there exists $\a\in\rtimes$ such that $(\a,\a)=1$. Then by
Proposition \ref{v2}($ii$) and $(T5)(a)$, we get
\begin{equation}\label{1.9}
(\gamma,\b)\in \bbbq\hbox{\;\;\;for all \;} \b\in
R,\gamma\in\rtimes.
\end{equation}

In the following proposition, we record several important
consequences of axioms $(T1)-(T5)$.
\begin{pro}\label{m1}  Let $L$ satisfy $(T1)-(T5)$. Then

(i) $(R, \rzero)=\{0\}.$

(ii) $(\a,\b)\in\bbbq$ for all $\a,\b\in R.$

(iii) For $\a\in R$, we have $t_\a\in\LLc.$ Moreover $t_\d\in
Z(\LLc)$ for all $\d\in R^0.$

(iv) For  $\a,\b,\gamma\in R^\times$ with $\b+\gamma\neq0$, we
have  $\hh_\a\subseteq \LLc$ and
$(\hh_\a,\LL_\b)=0=(\hh_\a,[\LL_\b,\LL_\gamma]).$

(v) For $\b\in R$ and $\a\in\rtimes$, there exist $d, u\in
{\mathbb Z}_{\geq 0}$ such that $d-u=2{(\b,\a)}/(\a,\a)$ and for
each $n\in\bbbn,$ $\b\tiny{\hbox{$+$}} n\a\in R$ if and only if
$-d\leq n\leq u.$

(vi) For $\b\in R$ and $\a\in\rtimes$, $-4\leq
\frac{2(\b,\a)}{(\a,\a)} \leq4.$

(vii) $\LLc$ is a perfect ideal of $\LL.$

(viii) $Z(\LLc)=rad(\cdot,\cdot)|_{\LLc\times \LLc}$ and
$\LLc^\bot=C_\LL(\LLc)$ where $\LLc^\bot:=\{x\in\LL\mid
(x,\LLc)=0\}.$

(ix) $Z(\LLc)\subseteq\bigoplus_{\d\in R^0}\LL_\d.$

(x) $\LL$ is tame if and only if $\LLc^\bot=Z(\LLc).$
 \end{pro}
\noindent{\bf Proof.} $(i)$ We first  prove $(\rtimes,
R^0)=\{0\}$. Suppose to the contrary that $\a\in\rtimes$, $\d\in
R^0$ and $(\a,\d)\not=0$. Since $(T3)$ holds, one can use the same
argument as in the  proof of   [AABGP, Lemma I.1.30] to show that
$\a+n\d\in R$ for infinitely many $n\in\bbbz$. But at most for one
$n$, $\a+n\d\in R^0$. So by   (\ref{1.9}) and  Proposition
\ref{v2}, $(\d,\a+n\d),(\a+n\d,\a+n\d)\in \bbbq$ and
$2(\d,\a+n\d)/(\a+n\d,\a+n\d)\in\bbbz $ for infinitely many
$n\in\bbbz$ which is impossible. Thus $(\rtimes , R^0)=\{0\}.$ Now
using this and ${( T5)}(b)$, it follows easily that
$(R^0,R^0)=\{0\}$.

$(ii)$ Immediately follows from  part $(i)$ and  (\ref{1.9}).

 $(iii)$ By   $(T3)$,  we have $t_\a\in\LL_c $ for
$\a\in\rtimes.$ Now  if $\d\in R^0$, then by $(T5)(b)$  and part
$(i)$ there exists $\a\in\rtimes$ such that $\a+\d\in\rtimes.$ So
we have $t_\a+t_\d=t_{\a+\d}\in\LL_c$. But $t_\a\in\LL_c$ and so
$t_\d\in\LL_c. $ The last statement follows from part $(i)$.

$(iv)$ Let $x\in\LL_\a$ and $y\in\LL_{-\a}$. Using  (\ref{m0}), we
have $[x,y]=(x,y)t_\a+h_{x,y}$ and so part $(iii)$ gives that
$h_{x,y}\in \LLc.$  Thus $\hh_\a=\hbox{span}_\bbbf\{h_{x,y}\mid
x\in\LL_\a,y\in\LL_{-\a}\}\subseteq \LLc.$ For the second
statement note that $h_{x,y}\in\LL_0$, so for  $b\in\LL_\b$,  we
have
\begin{eqnarray*}
(h_{x,y},b)=(h_{x,y},\frac{1}{2}[h_\b,b])=\frac{1}{2}([h_{x,y},h_\b],b)=
\frac{1}{2}(0,b)=0.
\end{eqnarray*}
Thus $(\hh_\a,\LL_\b)=0.$ Now suppose  $x\in\LL_\a$,
$y\in\LL_{-\a}$, $a\in\LL_\b$ and $b\in\LL_\gamma$. Using
(\ref{v4}), we have
$$(h_{x,y},[a,b])=([h_{x,y},a],b)\in(\LL_\b,\LL_\gamma)=0.$$ In other words
$(\hh_\a,[\LL_\b,\LL_\gamma]) =0$.

($v$) See [AKY, Proposition 1.7(a)].

($vi$) Since we have already  assumed $(\a,\a)=1$ for some
$\a\in\rtimes$, we have
 $(\gamma,\gamma)>0$ for all  $\gamma\in \rtimes$  [AKY, Proposition 1.8].
  Now use the same argument as in [AABGP, Lemma I.2.6].

($vii$) We know that  $\LL=\op_{\a\in R}\LL_\a=\sum_{\a\in\rtimes
}\LL_\a\op\sum_{\a\in\rzero}\LL_\a$. So we have
\begin{equation}\label{score}
\LLc=\sum_{\a\in\rtimes }\LL_\a\op\sum_{\hbox{\tiny
$\begin{array}{c}
\a,\b\in\rtimes\\
\a+\b\in\rzero \end{array}$}}[\LL_\a,\LL_\b].
\end{equation}
Now using part $(iii)$, we are done.

$(viii)$ One can see this,  using the same argument as in [BGK,
Lemma 3.6].

 $(ix)$ It follows from (\ref{score}) that
 \begin{equation} \LLc=\op_{\a\in R} (\LL_\a\cap\LLc)\andd
 Z(\LLc)=\op_{\a\in R}(Z(\LLc)\cap\LL_\a). \label{v3}
\end{equation} Now let $\a\in\rtimes$ and $x\in
Z(\LLc)\cap\LL_\a,$ then
part $(i)$ together with (\ref{v4}) gives   that $(x, \LL_{\b})=0$
for all $\b\in R$.
So the non-degeneracy of the form implies  that $x=0$.   Using
(\ref{v3}), we have
$Z(\LLc)\subseteq\bigoplus_{\d\in R^0}\LL_\d.$

($x$) This follows from   part ($viii$).
  \qed



\vspace{.3cm}

 Now for    $\LL\in\T$, denote by  $\v_\bbbq$   the $\bbbq-$span of the root system $R$ of
 $\LL$ and set
$\v:=\bbbr\ot_\bbbq\v_\bbbq.$ Using Proposition \ref{m1}$(ii)$,
one can get, in a natural way, an induced real valued form on $\v$
denoted by $(\cdot,\cdot)_\bbbr. $ Now identify $\v_\bbbq$ as a
subset of $\v$. In the appendix we prove the following theorem:
\begin{thm}
Let $\LL\in\T$ with root system $R$. Then $(R,
(\cdot,\cdot)_\bbbr)$ is an extended affine root system in $\v$ in
the sense of [AKY] (see the Appendix).\label{lr}
\end{thm}

\section{ The structure of  $\LLc$ } In 1992,  Berman and  Moody  [BM]
 introduced the class of $\D-$graded Lie algebras. They classified  the Lie algebras graded
by irreducible  simply-laced finite root systems  of rank $\geq2,$
i.e. root systems of type $A_\ell(\ell\geq2)$,
$D_\ell(\ell\geq4),$ $E_6,$ $E_7,$  $E_8$. In 1996,  Benkart and
Zelmanov [BZ] did the classification for the remaining reduced
types namely  for types $A_1,$ $B_\ell$, $C_\ell$, $F_4$ and
$G_2$. Finally,  Allison, Benkart and  Gao generalized the concept
of $\D-$graded Lie algebras to non-reduced types [ABG]. A Lie
algebra graded by an irreducible finite root system $\D$ consists
of a finite dimensional split simple Lie algebra which is of type
$\D$ if $\D$ is reduced and of type $B, C$ or $D$ if $\D$ is
non-reduced. In [Yo2],  Yoshii considered $\D-$graded Lie algebras
whose corresponding finite dimensional split simple Lie
subalgebras are of type $\D_{red}$. Yoshii's concept comes from
the theory of EALAs, namely, the core of an EALA is a $\D-$graded
Lie algebra in the sense of Yoshii [Az1, AG]. Yoshii also
introduced $(\D, G)-$graded Lie algebras, $\D$ an irreducible
finite root system and $G$ an abelian group [Yo2]. In this section
we prove that the core of an element in $\T$ is a $(\D,G)-$graded
Lie algebra for an irreducible  finite root system $\D$ and a free
abelian group $G$. Similar to [AABGP, Chapter III, \S1], we also
give a general construction of elements in $\T$ which  satisfy a
{\it division property},  starting from a Lie algebra $\gg$
satisfying certain properties. In fact we prove that if $\LL\in\T$
satisfies this division property, then $\gg$ is isomorphic to
$\LL_{cc}$. Throughout this section $\D$ is an irreducible finite
root system and $G$ is an abelian group.

\subsection{  $(\D,G)-$graded Lie algebras}
Let $\mathfrak{g}$ be a finite dimensional split simple Lie
algebra over $\bbbf$ with a splitting  Cartan subalgebra $\fh$ and
root system $\D_{red}$ so that $\fg$ has a root space
decomposition $\fg=\op_{\mu\in\D_{red}}\fg_{\mu}$ with
$\fh=\fg_0.$
\begin{DEF}
{\em Let $\fg$ and $\fh$ be as above. A $\D-$graded Lie algebra
$\LL$  over $\bbbf$ with grading pair $(\fg,\fh)$ is  a Lie
algebra satisfying the following conditions:

 (i) $\LL$ contains $\fg$ as a
subalgebra,

(ii) $\LL=\op_{\mu\in\D}\LL_{\mu}, \hbox{\;where\;}
\LL_{\mu}=\{x\in\LL\mid [h,x]=\mu(h)x\hbox{\;for all\;}h\in\fh\},$

(iii) $\LL_0=\sum_{\mu\in\D^\times}[\LL_{\mu},\LL_{-\mu}]. $}
\end{DEF}
\begin{DEF}{\em  A $\D-$graded Lie algebra
$\LL=\op_{\mu\in\D}\LL_{\mu}$ with grading pair $(\fg,\fh)$ is
called {\it  $(\D,G)-$graded }if $\LL=\op_{g\in G}\LL^g$ is a
$G-$graded Lie algebra such that $\fg\sub \LL^0$ and
supp$(\LL):=\{g\in G\mid \LL^g\neq\{0\}\}$ generates $G$. Since
$\fg\subseteq\LL^0,$ $\LL^g$ is an $\fh-$module for $g\in G$ and
by [MP, Proposition 2.1] we have $\LL=\op_{\mu\in\D}\op_{g\in
G}\LL^g_\mu$ where $\LL_\mu^g:=\LL^g\cap\LL_\mu$ for $g\in G$ and
$\mu\in\D.$}
\end{DEF}
\begin{DEF}
{\em Let $\LL=\op_{\mu\in\D}\op_{g\in G}\LL^g_\mu$ be a
$(\D,G)-$graded Lie algebra with grading pair $(\fg,\fh)$.  $\LL$
is called a {\it division $(\D,G)-$graded} Lie algebra if for each
$\mu\in\D^\times$, $g\in G$ and  $0\neq x\in\LL^g_\mu$, there
exists $y\in\LL^{-g}_{-\mu}$ such that $[x,y]=t_\mu (\mod\;
Z(\LL))$. A division $(\D,\bbbz^n)-$graded Lie algebra $\LL$ with
$\dim_\bbbf(\LL_\mu^\sg)\leq1$ for all $\sg\in\bbbz^n$ and
$\mu\in\D^\times$ is called a {\it Lie $n-$torus}  or  simply a
{\it Lie torus.}}
\end{DEF}
Now let $\LL\in \T$  with root system $R$ and nullity $\nu$, then
by Theorem \ref{lr}, $(R, (\cdot,\cdot)_\bbbr)$ is an extended
affine root system in $\v=\bbbr\ot_\bbbq\v_\bbbq$. Take $\v^0$ to
be the radical of the form $(\cdot,\cdot)_\bbbr$, $\bar{\v}=
\v/\v^0$ and $^-:\v\longrightarrow\bar{\v}$  to be the canonical
map. By [AKY, \S2], the image $\bar{R}$ of   $R$ under $^-$, is an
irreducible finite root system in $\bar{\v}$. We fix a choice of a
fundamental system $\{\bar{\a}_1,\ldots,\bar{\a}_\ell\}$ for
$\bar{R}.$ Also for each $i,$ we fix a preimage $\adot_i$ in $R$
of $\bar{\a}_i$ under $^-.$ Then there exists an irreducible
finite root system $\dot{R}$ with base
$\{\adot_1,\ldots,\adot_\ell\}$ in
$\dot{\v}:=\hbox{span}_\bbbr\{\adot_i\}_{i=1}^\ell$ such that $
\dot R$ is isomorphic to $\bar{R}$ and $\rd_{red}\subseteq R$.
Also $\v=\dot \v\op \v^0$ and $R^0= R\cap\v^0.$  The {\it rank}
and the {\it type} of $\LL$  is defined to be the rank and the
type of $\dot{R}$. Let $\rds$, $\rdl$ and $\rde$ be the set of
short, long and extra long roots of $\dot{R}$, respectively. Set
$$\begin{array}{c}S:=\{\d\in \v^0\mid\d+ \rds\sub R\},\quad
L:=\{\d\in\v^0\mid\d+ \rdl\sub R\},\hbox{ if }\rdl\not=\emptyset,\\
\hbox{ and}\;\; E:=\{\d\in \v^0\mid\d+ \rde\sub R\},\hbox{ if
}\rde\not=\emptyset. \end{array}$$ Then $S,L$ are semilattices in
$\v^0$ and $E$ is a translated semilattice in $\v^0$  (see [AABGP,
Chapter II] for the terminology). Choose a fixed subset
$\{\d_j\}_{j=1}^\nu$ of $S$ such that
$\vz=\sum_{j=1}^\nu\bbbr\d_j$  [AABGP, Proposition II.1.11]. Now
set $\dot{\hh}:=\sum_{i=1}^\ell\bbbf t_{\dot\a_i}$ and
$\hh^0:=\sum_{j=1}^\nu\bbbf t_{\d_j}$.
By   Proposition \ref{m1}$(iii)$, we have
\begin{equation}\label{y0}
\hh^0\subseteq Z(\LLc)\andd\dot{\hh}\op\hh^0\subseteq\LLc.
\end{equation}
If $\adot\in \rd_{red}\subset \rtimes$, then by  $(T3)$ there
exist $e_{\pm\adot}\in\LL_{\pm\adot}$ such that $(e_{\adot},
h_{\adot}, e_{-\adot})$ is an $\mathfrak{sl}_2\hbox{-triple}.$
Also we know that for  $1\leq i\leq\ell,\; \adot_i\in
\dot{R}_{red}\sub \rtimes$, so we can define
\begin{equation}\label{m11}
\dg:=\hbox{subalgebra of\;}\LLc\hbox{\;generated by \;}
\{e_{\pm\adot_i}\}_{i=1}^\ell.\end{equation} For $1\leq
i,j\leq\ell$, put $c_{i,j}:=\adot_j(h_{{{\adot}_i}})$. Since by
axiom $(R4)$ (see  Appendix, Definition \ref{ears}),
$\adot_{j}+(-c_{i,j}+1)\adot_i\not\in R$ for $1\leq i,j\leq\ell,$
the following relations are satisfied in $\dg:$
\begin{equation}\label{m13} \begin{array}{c}
[h_{{{\adot}_i}}, h_{{{\adot}_j}}]=0,\; [h_{{{\adot}_i}},
e_{{{\adot}_j}}]=c_{i,j}e_{{{\adot}_i}},\;
[h_{{{\adot}_i}},e_{{{-\adot}_j}}]=-c_{i,j}e_{{{-\adot}_i}},\;
 \\
\ [e_{{{\adot}_i}}, e_{{{-\adot}_j}}] =\d_{i,j}h_{{{\adot}_i}},\;
\theta^+_{i,j}=0,\;
 \theta^-_{i,j}=0,\\
\end{array}
1\leq i,j\leq\ell
\end{equation}
where $\theta^+_{i,j}=(\ad
e_{{{\adot}_i}})^{-c_{i,j}+1}(e_{{{\adot}_j}})\hbox{\; and\;}
\theta^-_{i,j}=(\ad
e_{{{-\adot}_i}})^{-c_{i,j}+1}(e_{{{-\adot}_j}}).$ Set
$C:=(c_{i,j})_{1\leq i,j\leq \ell}$, then $C$ is the Cartan matrix
corresponding to $\rd_{ red}$ and by  Serre's Theorem (see e.g.
[MP, Proposition 4.3.3, Theorem 4.6.4]) $\dg$ is a  finite
dimensional split simple Lie algebra. Moreover
$\{h_{{\adot}_i}\}_{i=1}^\ell$ is a linearly independent subset of
$\hh$ which forms a basis for  a splitting  Cartan subalgebra of
$\dg$. Therefore we have the following theorem:
\begin{thm}\label{m12}
$\dg$ is a finite dimensional split simple subalgebra of $\LLc$
with root system $\rd_{ red}$ and splitting  Cartan subalgebra
$\dhh=\op_{i=1}^\ell\bbbf h_{\adot_i}$.
\end{thm}
With the same  notations as before, the following lemma holds:

\begin{lem}
Let $\LL\in\T,$ then the form $(\cdot,\cdot)$ restricted to
$\dot{\hh}$ is non-degenerate. \label{ac}
\end{lem}
\noindent{\bf Proof.}
Let $\adot\in\rd_{red}.$ By Theorem \ref{m12},
$\dg_{\adot}\subseteq\LL_{\adot}.$ Also, it is known from the
finite dimensional theory that $0\neq
[\dgg_{\adot},\dgg_{-\adot}]\subseteq \dot{\hh}$. So (\ref{m0})
together with the fact that
$\dim(\dgg_{\adot})=\dim(\dgg_{-\adot})=1$ implies that there
exist $x\in\dgg_{\adot}$ and $y\in\dgg_{-\adot}$ so that $0\neq
[x,y]=(x,y)t_{\adot}$. This means that the form $(\cdot,\cdot)$
restricted to $\dgg$ is nontrivial and so  is non-degenerate [MP,
Exercise 1.7]. Also since $\dg$ is a finite dimensional split
simple Lie algebra, it is central simple. So the form on $\dg $ is
a non-zero scalar multiple of the Killing form.
Therefore $(\cdot,\cdot)|_{\dhh\times\dhh}$ is non-degenerate.\qed

\begin{rem}\label{rema}
Since $t_{\adot}\in\dot\hh=\sum_{i=1}^\ell\bbbf t_{\adot_i}$ for
all $\adot\in\rd\subseteq\hh^\star,$ the unique element in
$\dot\hh$ representing $\adot$ through the form  restricted to
$\dot\hh$ is the same as  the unique element in $\hh$ representing
$\adot$ through the form on $\hh.$
\end{rem}
Since  $ \dim(\hh)<\infty,$ $\hh^0\subseteq\hh$ and
$(\dot{\hh}\op\hh^0, \hh^0)=0,$ there exists a subspace $D$ of
$\hh$ such that
\begin{equation}\label{t}
\begin{array}{l}
\dim D=\dim \hh^0, \\
(\dot{\hh}\op D, D)=0\hbox{ and,}\\
\hbox{  the form is non-degenerate on\;} \dot{\hh}\op\hh^0\op D.
\end{array}
\end{equation}Let   $W$ be the orthogonal complement on $\dot{\hh}\op\hh^0\op D$ in $\hh$ with respect to the form, then
\begin{equation}
\hh=\dot{\hh}\op\hh^0\op D\op W\andd\LL_0=\dot{\hh}\op\hh^0\op
D\op W\op\hh^\bot. \label{G2}
\end{equation}
\begin{lem}\label{m3} With the above notations we have

(i) $W\subseteq Z(\LL)\subseteq C_\LL(\LLc).$

(ii) $\LLc\cap\hh=\dot{\hh}\op\hh^0\hbox{\;and\;}
Z(\LLc)\cap\hh=\hh^0.$

 (iii) $ \hbox{If\;} \LL \hbox{\;is tame,
then\;} W=0. $
\end{lem}
\noindent{\bf Proof.} $(i)$ Let $h\in W$ and  $\a\in R,$ then
$(t_\a, h)\in(\dot{\hh}\op\hh^0, W)=0.$ So if $x\in\LLa,$ then
$[h,x]=\a(h)x=(t_\a,h)x=0$ which implies $W\subseteq Z(\LL).$

$(ii)$ Let $d\in D, w\in W$ and $d+w\in\LLc.$ By $(i)$ and
Proposition \ref{m1}($viii$), $( W,w)=( W,w+d)=\{0\}$, so $w\in
\hh^\bot\cap W=\{0\}$. This together with (\ref{m0}) and
(\ref{G2}) implies that  $d=0$. Therefore $(D+W)\cap\LLc=\{0\}$
and then (\ref{y0}) implies that $\LLc\cap\hh=\dot{\hh}\op\hh^0.$
For the second statement,  it is enough to show $\dot{\hh}\cap
Z(\LLc)=\{0\}$ as by (\ref{y0}), $\hh^0\subseteq Z(\LLc)$. If
$x\in \dot{\hh}\cap Z(\LLc),$ then Proposition \ref{m1}($viii$)
implies that $(x, \LLc)=\{0\},$ hence we have  $(x,
\dot{\hh})=\{0\}.$ But the form is non-degenerate on $\dot{\hh}$,
so $x=0.$

$(iii)$ If $\LL$ is tame, then by part $(i)$ and Proposition
\ref{m1}($viii$), $W\subseteq C_\LL(\LLc)=\LLc^\bot=Z(\LLc)$.
Therefore $W\sub\LLc\cap\LLc^\bot$ and so  $(W,W)=0.$ But by
(\ref{G2}), the form on $W$ is non-degenerate. Thus $W=0.$\qed
\\
$ $

Consider (\ref{G2}), $\sdhh$ is imbedded in $\hhstar$ by extending
 a linear form  $\a\in \sdhh$ by zero on the other three summands, and similarly  for
 ${(\hhz)}^{\star},\;  D^\star,\; W^\star.$ So we can identify
$$\hhstar=\sdhh\op{(\hhz)}^{\star}\op D^\star\op W^\star.$$
Contemplating (\ref{v3}), one can see that  $\dhh\subseteq\hh$ is
an abelian subalgebra of $\LL$ which is ad-diagonalizable on
$\LLc$. So we have
$\LLc=\sum_{\a\in\hhstar}(\LLc\cap\LL_\a)=\sum_{\adot\in\dot{\hh}^\star}(\LLc)_{\dot{\a}}$,
where
\begin{equation}\label{b00}
\begin{array}{c}
(\LLc)_{\dot{\a}}=\{x\in\LLc\mid
[h,x]=\adot(h)x,\;h\in\dhh\};\;\;\adot\in\dhh^\star.
\end{array}
\end{equation}

Now let $\pi:\hhstar\rightarrow\sdhh$ be the projection map. We
have the  following lemma:
\begin{lem}\label{m4}
$ \{\adot\in\sdhh\mid(\LLc)_{\adot}\neq\{0\}\}=\pi(R)=\dot{R}.$
\end{lem}

\noindent {\bf Proof.} Let $\gamma \in R$. Then
$\gamma=\dot\gamma+\d$  for some $\dot\gamma\in\rd$ and $\d\in
R^0$ (see Appendix, (\ref{q4}) and Lemma \ref{aa}). Then for
$h\in\dot\hh$, we have
$$(t_{\dot\gamma},h)=(t_{\dot\gamma}+t_\d,h)=(t_\gamma,h)=\gamma(h)=\pi(\gamma)(h).$$But
the form on $\dot\hh$ is non-degenerate (Lemma \ref{ac}), so
$t_{\dot\gamma}=0$ and consequently  $\gamma=\dot\gamma+\d=\d\in
R^0.$ It means $\pi(\gamma)=0$ if and only if $\gamma\in R^0.$ Now
we are done using this together with the fact that
$(\LLc)_{\dot{\b}}=\sum_{\{ \a\in R;\pi(\a)=\dot{\b}\}}(\LLc\cap
\LL_\a)$ for $\dot\b\in \dot\hh^\star.$  \qed
\\

 We recall that $\Lam=\la R^0\ra$. Then  $\LL$ is
a $\Lam-$graded Lie algebra, namely
\begin{equation}\label{m8}
\LL=\op_{\sg\in\Lam}\LL^\sg\hbox{\;\;where\;\;}
\LL^\sg=\sum_{\adot\in\dot{R}}\LL_{\adot+\sg},\;\sg\in\Lam.
\end{equation}
Since $\LLc$ is an ideal of $\LL$ generated by homogeneous
elements with respect to this grading, $\LLc$ is a $\Lam-$graded
ideal and so
\begin{equation}\label{gradcore} \LLc=\op_{\sg\in\Lam}\LLc^\sg
\hbox{\;\;where\;\;} \LLc^\sg=\LL^\sg\cap\LLc,\;\sg\in\Lam.
\end{equation}
For $\adot\in\rd$ and $\tau\in\Lam$ set
$(\LLc)_{\adot}^\tau:=(\LLc)_{\adot}\cap\LLc^\tau.$
It is easy to see that
\begin{equation} \label{m9}
(\LLc)_{\adot}=\sum_{\tau\in\Lam}(\LLc)_{\adot}^\tau \andd
(\LLc)_{\dot\b}^\sg=\LL_{\dot\b+\sg}
;\;\;\adot\in\rd,\;\dot\b\in\rd^\times,\sg\in\Lam.
\end{equation}

\begin{thm}\label{deltagraded}
Let $\LL\in\T.$ Then $\LLc$ is a division  $(\rd,\Lam)-$graded Lie
algebra with grading pair $(\dg,\dot\hh)$ where $\dgg$ is defined
as  in (\ref{m11}).
\end{thm}
\noindent{\bf Proof.}  By Theorem \ref{m12}, $\dgg$ is a finite
dimensional split simple subalgebra of $\LLc$ with splitting
Cartan subalgebra $\dhh$. Since
$e_{\pm\adot_i}\in\LL_{\pm\adot_i}$ for $1\leq i\leq\ell$,
$\dgg_{\adot}\subseteq\LL_{\adot}$ for all $\adot\in\rd_{red}$. So
using (\ref{m8}) and (\ref{gradcore}), we have
\begin{equation}\label{m14}
\dgg=\op_{\adot\in\dot{R}_{red}}\dgg_{\adot}\subseteq(\op_{\adot\in\rd_{
red}}\LL_{\adot+0})\cap\LLc\subseteq\LLc^0.
\end{equation}
Now (\ref{b00}) together with Lemma \ref{m4} implies that
$$
\LLc=\op_{\adot\in\rd}(\LLc)_{\dot{\a}}\;\hbox{where}\;(\LLc)_{\dot{\a}}=\{x\in\LLc\mid
[h,x]=\adot(h)x,\;h\in\dhh\}. $$ Also  using (\ref{score}), one
can easily see that
\begin{equation}(\LLc)_{0}=\sum_{\dot{\eta}\in\rd^\times }[(\LLc)_{\dot{\eta}},
(\LLc)_{-\dot{\eta}}]. \label{m16}
\end{equation} Therefore $\LLc$ is an $\rd-$graded Lie algebra.
On  the other hand  $\LLc$ is a $\Lam-$graded Lie algebra with
$\dg\subset\LLc^0$  and  such that supp$(\LLc)=R^0$. Thus $\LLc$
is an $(\dot{R},\Lam)-$graded Lie algebra. Finally   the division
property is a consequence of (\ref{m9}), $(T3)$ and Proposition
\ref{m1}$(iii).$ \qed
  \begin{cor}
 Let $\LL\in\T.$ Then $\LL_{cc}$ is a centerless  division $(\rd,\Lam)-$ graded Lie algebra with
grading pair $(\dg,\dot\hh).$ In particular $\LL_{cc}$ is a
$\Lam-$graded simple Lie algebra.
 \label{c02}
 \end{cor}
 \noindent{\bf Proof.}  The natural $\dot {R}$ and $\Lam-$grading on $\LL_{cc}$ induced
 from  $\LLc$  together with  Theorem \ref{deltagraded} show that
 $\LL_{cc}$ is a centerless  division $(\rd,\Lam)-$ graded Lie
 algebra.
 The second statement  is a consequence of  [Yo2, Lemma
 4.4].\qed

\subsection{ A characterization  for  a subclass of $\T$}  Let $\LL\in\T$ and $^-:\LLc\rightarrow\LL_{cc}$ be the canonical map. By
Proposition  \ref{m1}$(ix)$ and Lemma \ref{m3}$(ii)$, $^-|_{\dgg}$
is injective, so we may identify $\dgg$ as a subalgebra of
$\LL_{cc}$. In this case we have the following lemma which is a
straightforward consequence of Proposition \ref{m1}$(iii)$:

\begin{lem}\label{m17}
 If $\adot\in\rd$ and  $\sg\in R^0$ with $\a=\adot+\sg\in R$, then  $t_{\adot}=\overline{t_\a}$.
\end{lem}

\begin{DEF}
{\em Let $\nu$ be a positive integer and let $\dd$ be the class of
all triples $(\gg,(\cdot,\cdot),\dot\fh)$, $\gg$  a Lie algebra,
$(\cdot,\cdot)$ a symmetric  bilinear form on $\gg$ and $\dot\fh$
a subalgebra of $\gg$, satisfying the following 12 axioms:

\noindent{\bf (D1)} The form is non-degenerate and invariant.

\noindent{\bf (D2)} $\dot\fh$ is a  finite dimensional  toral
subalgebra of  $\gg$. Take $\rd$ to be the root system of $\gg$
with respect to $\dot\fh.$

 \noindent{\bf (D3)} $(\cdot, \cdot)|_{\dot\fh\times\dot\fh}$ is non-degenerate and the form  $(\cdot,\cdot)$ restricted to the $\bbbq-$ subspace of
 $\dot\fh$ spanned by  $\{t_{\adot}\mid \adot\in\rd\}$ is $\bbbq-$valued.

 ${ (D3)}$ allows us to transfer  the form  on $\dot\fh$ to a
form on
 $\dot\fh^\star$ by setting  $(\adot,\dot{\b}):=(t_{\adot},t_{\dot{\b}})$ for
$\adot,\dot{\b}\in\dot\fh^\star$. Put
$\ud_\mathbb{Q}:=\hbox{span}_{\mathbb{Q}}(\dot R)$ and
$\ud:=\mathbb{R}\otimes_\bbbq\ud_\mathbb{Q}.$ Consider a basis
$\{\dot{\b}_i\}_{i\in I}$ of $\ud_\bbbq$ and define
$(1\ot\dot{\b}_i,1\ot\dot{\b}_j):=(\dot{\b}_i,\dot{\b}_j)$ for
$i,j\in I$. Extend this linearly to a bilinear form on $\dot{\u}$.
We may consider $\rd$ as a subset of $\dot{\u}$ by identifying
$\adot$ with $1\ot\adot.$ The next axioms are:

\noindent{\bf (D4)} The  form on $\ud$ is a positive definite
symmetric bilinear form and  $\dot R$ is an irreducible finite
root system in $\ud.$

\noindent{\bf (D5)} $\gg=\op_{\sg\in\zn}\gg^\sg$ is a $\zn-$graded
Lie algebra.

\noindent{\bf (D6)} For each $\adot\in\dot R$,
$\gg_{\adot}=\op_{\sg\in\zn}\gg^\sg_{\adot}$ where
$\gg^\sg_{\adot}:=\gg^\sg\cap\gg_{\adot}$ for $\sg\in\zn$.

\noindent{\bf (D7)} $\dot\fh\subseteq\gg^0\cap\gg_0=\gg_0^0.$

\noindent{\bf (D8)}
$\gg_0=\sum_{\adot\in\dot{R}^\times}[\gg_{\adot}, \gg_{-\adot}]$.

\noindent{\bf (D9)} $\{\sg\in\zn\mid\gg^\sg\neq\{0\}\}$ generates
a subgroup of $\zn$ of rank $\nu$.

\noindent{\bf (D10)}  If  $\sg,\tau\in\hbox{supp}(\gg)$, then
$(\gg^\sg, \gg^\tau)=\{0\}$ unless $\sg+\tau=0.$

Note that  for $\dot\a,\dot\b \in\rd$,  the invariance of the form
and the Jacobi identity imply
 that
$(\gg_{\dot\a},\gg_{\dot\b})=\{0\}$ unless $\dot\a+\dot\b=0$. Now
using this together with $({D10})$, we have
$(\gg_{\dot\a}^\sg,\gg_{\dot\b}^\tau)=\{0\}$ if $\sg+\tau\neq0$ or
$\dot\a+\dot\b\neq0$.

\noindent{\bf (D11)} $\gg^0\cap\gg_{\adot}\neq\{ 0\}$, for all
$\adot\in \rd_{red}^\times.$

\noindent{\bf (D12)} ($a$) For each $\adot\in \dot{R}^\times$,
$\sg\in\zn$ and $0\neq x\in \gg^\sg\cap\gg_{\adot}$ there exists
$y\in\gg^{-\sg}\cap\gg_{-\adot},$ such that $[x, y]=t_{\adot}$.

($b$) For $\sg\in\zn$ if  $ \gg^\sg\cap\gg_{0}\neq \{0\}$ there
exist $x\in \gg^\sg\cap\gg_{0}$ and $y\in\gg^{-\sg}\cap\gg_{0}$
such that $[x, y]=0$ and $(x,y)=1$.

When $(\cdot,\cdot)$ and $\dot\fh$ are fixed or unimportant, we
denote a triple $(\gg,(\cdot,\cdot),\dot\fh)$ by $\gg$. }
\end{DEF}

The property stated in $(D10)$ says that the form on $\gg$ is
$\zn-$graded in the following sense:
\begin{DEF}
{\em  A  symmetric bilinear form $(\cdot,\cdot)$ on a $G-$graded
Lie algebra $\LL=\op_{g\in G}\LL^g$ satisfying
$$(\LL^g,\LL^h)=0\;\hbox{unless}\; g+h=0\; \hbox{for all}\;g,h\in
\hbox{supp}(\LL),$$is called a  {\it $G-$graded form} (or a {\it
graded form} for simplicity).}
\end{DEF}

Let    $\LL\in\T$ and  put $\LL_{cc}:=\LLc/Z(\LLc)$ Consider the
split simple Lie algebra $\dot\gg$ with splitting Cartan
subalgebra $\dhh$ as in Theorem \ref{m12}. We recall that the
canonical map $^-:\LLc \longrightarrow\LL_{cc}$ restricted to
$\dot\gg$ is injective and so we identify $\dot\gg$ as a
subalgebra of $\LL_{cc}.$ It follows from Proposition
\ref{m1}($viii$) that the form on $\LL$ induces a form
$(\cdot,\cdot)_{cc}$ on  $\LL_{cc}$ such that
$(\bar{x},\bar{y})_{cc}=(x,y)$ for $x,y\in\LLc$. In the following
theorem, we give a characterization  of algebras in $\T$ in terms
of  algebras in $\dd$. Our argument is essentially the same as
[Az1, Proposition 1.28], however as our axioms  are different from
those of [Az1], for the convenience of reader, we provide the
proof.
\begin{thm}\label{main}
Let $(\LL,(\cdot,\cdot),\hh)\in\T$  with root system $R$ and
nullity $\nu$. Then $(\LL_{cc},(\cdot,\cdot)_{cc},\dot{\hh})$
satisfies $(D1)-(D11)$. Moreover  if $\LL\in\T$ satisfies the
following property
\begin{equation}\label{newp}
\begin{array}{c}
\hbox{if
 $\sg\in R^0$ with  $(\LLc)_{0}^\sg\neq\{0\}$, then there exist $x\in(\LLc)_{0}^\sg$ and
 $y\in(\LLc)_0^{-\sg}$}\\
 \hbox{ so that  $\overline{[x,y]}=0$ and
 $(x,y)=1$}
 \end{array}
\end{equation} then $({D12})$ is also satisfied. In other words, the
centerless core of a Lie algebra in $\T$ satisfying (\ref{newp})
  belongs to $\dd$.
\end{thm}
\noindent{\bf Proof.} We know that $\la R^0\ra$ is a free abelian
group of  rank $\nu$, so we may assume $\la R^0\ra=\zn.$ If
$x+Z(\LLc)\in\LL_{cc}$ and $(x+Z(\LLc), \LL_{cc})_{cc}=\{0 \}$,
then $(x, \LLc)=\{0\}$, so by Proposition \ref{m1}($viii$), we
have $x\in Z(\LLc).$ It means that $(., .)_{cc}$ is
non-degenerate. Also since $(\cdot,\cdot)$ is invariant,
$(\cdot,\cdot)_{cc}$ is invariant. So ({D1}) holds. By Corollary
\ref{c02}, $\LL_{cc}$ is an $(\rd,\zn)-$graded Lie algebra with
grading pair $(\dg,\dot\hh)$, therefore we have $({D2})$,
$({D5})$, $({D6})$, $({D8})$ and $(D11)$. From Proposition
\ref{m1}($ii$) we know that the form $(\cdot,\cdot)$ restricted to
$\v_\bbbq=\hbox{span}_{\bbbq}(R)$ is $\bbbq-$valued. Also from
subsection 2.1 we know $\bbbr\ot_\bbbr\v_\bbbq=\v=\dot\v\op\v^0$
where $\v^0$ is the radical of the induced form on $\v$ and
$\dot\v=\sum_{i=1}^\ell\bbbr{\adot}_i$ is a finite dimensional
real vector space  where $\{\adot_i\}_{i=1}^\ell$  is a base for
the finite  root system  $\dot{R} $ in $\dot\v.$  Also by Lemma
\ref{ac}, the form restricted to $\dot\hh$ is non-degenerate.
Therefore Remark \ref{rema} gives $(D3)$. Set
$\dot{\u}_\bbbq:=\hbox{span}_{\bbbq}{(\rd)}$ and
$\dot{\u}:=\bbbr\otimes_\bbbq\dot{\u}_\bbbq.$ We can identify
$\dot{\u}\equiv\dot{\v}$ and so $\rd$ is a finite root system in
$\dot{\u}$. This
 implies that $({D4})$ holds.



We know that $\{\adot_i\}_{i=1}^\ell$ is a subset of $R$. So by
$(T3)$ for each $1\leq i\leq\ell, $ there exist
$x_i\in\LL_{\adot_i}$ and $y_i\in\LL_{-\adot_i}$ such that
$t_{\adot_i}=[x_i, y_i],$ therefore $[x_i,
y_i]\in(\LLc)^0\cap(\LLc)_{0}$. Thus
$t_{\adot_i}=\overline{t_{\adot_i}}\in(\LL_{cc})_{0}^0$ (see Lemma
\ref{m17}) and so $\dhh\subseteq (\LL_{cc})_{0}^0.$ In other words
$(D7)$ holds.

Let $\adot\in\rds$ and  define $S=S_{\adot}:=\{\sg\in\rzero\mid
\adot+\sg\in R\}$. If $\sg\in S$, then $\adot+\sg\in\rtimes,$ so
$\LL_{\adot+\sg}\in\LL^\sg\cap\LLc.$ Therefore by Proposition
\ref{m1}($ix$),  $({\LL_{\adot+\sg}+Z(\LLc)})/{Z(\LLc)} $ is a
nonzero subset of $(\LL_{cc})^\sg$, hence
$S\subseteq\{\sg\in\zn\mid (\LL_{cc})^\sg\neq\{0\}\}$. Now since
$\la S\ra=\la R^0\ra=\zn,$ the set $\{\sg\in\zn\mid
(\LL_{cc})^\sg\neq\{0\}\}$ generates a subgroup of $\bbbz^\nu$ of
rank $\nu.$


Next  let $\sg,\tau\in R^0$ and  $\sg+\tau\neq0, $ then by
(\ref{v4})
$$((\LL_{cc})^\sg,(\LL_{cc})^\tau)_{cc}=(\LLc^\sg,\LLc^\tau)\subseteq(\sum_{\adot\in\rd^\times}\LL_{\adot+\sg}+\LL_\sg
,\sum_{\adot\in\rd^\times}\LL_{\adot+\tau}+\LL_\tau)=\{0\}.$$ It
means that $\LL_{cc}$ satisfies ${(D10)}$.

For the last statement it is enough to show that if $\sg\in\zn,$
$\adot\in\rd^\times$ and $0\neq
x+Z(\LLc)\in(\LL_{cc})^\sg_{\adot},$ then there exists
$y+Z(\LLc)\in (\LL_{cc})^{-\sg}_{-\adot}$ such that
$[x+Z(\LLc),y+Z(\LLc)]=t_{\adot}$.  Using (\ref{m9}), we have
$\LL_{\adot+\sg}\neq \{0\}$ which means $\a=\adot+\sg\in R$.
Without lose of generality we may assume $0\neq x\in\LL_{\a}$. So
$(T3)$ implies that there exists
 $y\in\LL_{-\a}$ such that $[x,y]=t_\a$. Then it
follows from Lemma \ref{m17} that
$[x+Z(\LLc),y+Z(\LLc)]=\overline{[x,y]}=\overline{t_\a}=t_{\adot}.$
 \qed
\begin{rem}
 Let $\LL\in\T$. Let $\tau,\lam\in\zn$ and $\dot \eta,\dot\b\in\rd$ so that
 $(\LL_{cc})_{\dot\eta}^\lam,(\LL_{cc})_{\dot\b}^\tau$ $\neq\{0\}.$ Theorem \ref{main} implies that the form $(\cdot,\cdot)_{cc}$ restricted to
 $(\LL_{cc})_{\dot\eta}^\lam\op(\LL_{cc})_{-\dot{\eta}}^{-\lam}$ is non-degenerate and
  $((\LL_{cc})_{\dot\eta}^\lam,(\LL_{cc})_{\dot\b}^\tau)_{cc}=0$
  except when $\dot{\eta}+\dot{\b}=0$ and $\lam+\tau=0$.
\end{rem}

\begin{pro}
Let $\LL=\op_{\mu\in \D}\op_{g\in G}\LL_{\mu}^{g}$ be a division
$(\D,\zn)-$graded Lie algebra with grading pair $(\fg,\fh)$ over
an algebraically closed field $\bbbf $ of characteristic zero.
Suppose that there is  a non-degenerate symmetric invariant graded
bilinear form on $\LL$  satisfying the following conditions:

$(i)$ The form restricted  to $\fg$ coincides with  the Killing
form.

$(ii)$ For  $\sg\in\zn$ and $\LL_0^\sg\neq0$, there exist
$x\in\LL_0^\sg$ and $ y\in\LL_0^{-\sg}$ such that $[x,y]=0$ and
$(x,y)=1. $

 Then $(\LL,(\cdot,\cdot),\fh)\in\dd$.\label{d}
\end{pro}
\noindent{\bf Proof.} The form is non-degenerate, symmetric and
invariant, so $({D1})$ holds. We have $({D2})$ because $\LL$ is a
$(\D,\zn)-$graded Lie algebra with grading pair $(\fg,\fh)$, in
fact $\fh$ plays the role of $\dot{\hh}$ in $({D2})$. It follows
using [H, \S 8.5] that the restriction of the Killing form to the
$\bbbq-$subspace $\fh$ spanned by $\{t_\mu\mid \mu\in\D_{red}\}$
is $\bbbq-$valued. So $(D3)$ holds. It is obvious that  $({D4})$
holds. $({D5}),({D6})$, $({D8})$ and $(D11)$  hold by the
definition of $\LL$  as a $(\D,\zn)-$graded Lie algebra. By [Yo3,
Lemma 2.1], $\LL$ is centerless, so if $\mu\in\D^\times$ and $g\in
G$ such that $\LL_\mu^g\neq 0$, the division property of $\LL$
implies that for $0\neq x\in \LL_\mu^g $ there exists $y\in
\LL_{-\mu}^{-g}$ such that $[x,y]=t_{\mu}$. This together with
$(ii)$ gives $(D12)$. We know that $\fg\subset\LL^0$, so
$\fh\subset\LL^0\cap\LL_0$ i.e. $({D7})$ is satisfied. $({D9})$
holds because supp$(\LL)$ generates $\zn$. Finally $\gg$ satisfies
$({D10})$ since the form is graded. \qed
\smallskip

Next we show how to construct elements in $\T$ of arbitrary
nullity starting from elements in $\dd$.  Consider $\gg\in\dd$ and
for
 $1\leq i\leq\nu$ define
 $d_i\in Der(\gg)$ such that
$d_ix=n_ix$ if $x\in\gg^{(n_1,\cdots, n_\nu)}.$ It follows  using
$({D9})$ that $d_1,\ldots,d_\nu$ are linearly independent over
$\bbbf.$

Set $\LL=\gg\op C\op D$ where $C=\op_{i=1}^\nu\bbbf c_i$ is a
$\nu-$dimensional vector space and $D=\op_{i=1}^\nu\bbbf
d_i\subseteq Der(\gg).$  Similar to [AABGP, $\S III.1$] we define
an anti-commutative product $[., .]'$  on $\LL$ and  also a
symmetric bilinear form $(. ,.)$  on $\LL$ which extends the form
on $\gg$ as follows:
\begin{equation}
\begin{array}{l}
\ [\LL, C]'=0, \\
\ [D, D]'=0,\\
\ [d_i, x]'=d_ix \;\;\; x\in\gg,\;1\leq i\leq\nu, \\
\ [x, y]'=[x, y]+\sum_{i=1}^\nu(d_ix, y)c_i\;\;\; x,y\in\gg,
\end{array}\label{ber}
\end{equation}
and
\begin{equation}
\begin{array}{l}
 (C, C)=0 , \;(D, D)=0,\\
(c_i, d_j)=\d_{i,j}, \;\;\;1\leq i,j\leq\nu,\\
(C, \gg)=(D, \gg)=0.
\end{array}\label{formm}
\end{equation}

The same argument as in [AABGP, Proposition III.1.20], with slight
modification, gives the following  theorem. We give a sketch  of
the proof and refer the interested  reader to  [AABGP] for
details.
\begin{thm}\label{cnonstruct}
Let $(\gg,(\cdot,\cdot),\dot\fh)$ be an element of $\dd$. Set
$\hh=\dot\fh\op C\op D$ and $\LL=\gg\op C\op D,$ then $(\LL, (.,
.), \hh)\in\T.$ Also $\LL$ is tame  of nullity $\nu$ and of the
same type as  the type of $\rd$. Moreover $\LLc=\gg\op C,$
$Z(\LLc)=C$ and $\LL_{cc}\simeq\gg$.
\end{thm}
\noindent{\bf Proof.} The form defined in (\ref{formm}) is a
non-degenerate  symmetric invariant bilinear form such that the
restricted form to $\hh$ is non-degenerate. It means that $(T1)$
and the last part of $(T2)$ hold. Let $\{\d_1,\ldots,\d_\nu\}$ be
the dual  basis  of $\{d_1,\ldots,d_\nu\}$ and identify
$\hh^\star$ with $\dot{\fh}^\star\op C^\star\op D^\star$. Then
$t_{\d_i}=c_i$ for  $i=1,\ldots,\nu$. We can identify
$\zn\subseteq D^\star$ by $(n_1,\ldots,n_\nu)=\sum_{i=1}^\nu
n_i\d_i.$ Then $\LL=\op_{\adot\in\rd,\sg\in\zn}\LL_{\adot+\sg}$
where
\begin{eqnarray}
\LL_{\adot+\sg} &=&\{x\in\LL\mid[h,x]'=(\adot+\sg)(h)x,\;h\in\hh\}\nonumber\\
&=&\left\{\begin{array}{ll}
(\gg_0\cap\gg^0)+C+D&\hbox{if}\;\; \adot+\sg=0\\
\gg_{\adot}\cap\gg^\sg&\hbox{otherwise.}
\end{array}\right.\label{rs}
\end{eqnarray}
Thus $\LL$ satisfies $(T2).$ The root system $R$ of $\LL$ is
$\{\adot+\sg\mid\adot\in\rd,\sg\in\zn,\gg_{\adot}\cap\gg^\sg\neq\{0\}\}\subseteq\rd+\zn.$
Moreover  $R^0= R\cap\zn$. Now ${(D12)}$ together with (\ref{m0})
and (\ref{ber}) implies that ${(T3)}$ holds. Since $\rd$ is a
finite root system, ${(T4)}$ holds. Set
$\u_\bbbq:=\hbox{span}_\bbbq(R)$ and
$\u:=\bbbr\ot_{\bbbq}\u_{\bbbq},$ then ${(D3)}$ allows us to have
a real valued form  $(\cdot,\cdot)$ on $\u$. If $\u^0$ is the
radical of the form $(\cdot,\cdot)$, then $\u=\dot\u\op\u^0$ where
$\dot\u=\bbbr\ot_\bbbq\hbox{span}_{\bbbq}(\rd)$. Now the canonical
map $^-:\u\longrightarrow\bar{\u}=\u/\u^0$ maps $\dot\u$
isometrically onto $\bar{\u}$. So by ${(D4)}$, $\bar{R}$, the
image
 $R$ under $^-,$ is  an  irreducible finite root system in $\bar{\u}$ which is isomorphic to $\rd.$ Therefore
 ${(T5)}(a)$ holds. Also ${(D8)}$ implies that ${(T5)}(b)$ holds. Finally $R^0\subseteq\zn,$ so $\la R^0\ra$
  is a free abelian group of finite rank. This completes the proof.\qed

\smallskip

We now provide an example of an element
$(\LL,(\cdot,\cdot),\hh)\in \T$ of type $C$ with arbitrary nullity
which satisfies all the axioms of an EALA except that $\hh$ is not
self centralizing. It is interesting that even though $\LL$  is
not  an EALA,  its root system remains an extended affine root
system in the sense of [AABGP]. In [AKY, Example 1.10], another
such an example, of type $BC,$ is provided but the root system is
not anymore an extended affine root system in the sense of
[AABGP].

 Let $(A,\cdot)$ be an algebra. A linear transformation $f:A\longrightarrow A$ is called an {\it
involution}, if the following conditions are satisfied
\begin{equation} f^2=id \andd
f(a\cdot b)=f(b)\cdot f(a) \;\hbox{for}\;a,b\in A. \label{in}
\end{equation}
If $V$ is a vector space over $\bbbc$, $f:V\longrightarrow V$ is called {\it semi-linear} if
$$
\begin{array}{ll}
\begin{array}{c}
f(v+w)=f(v)+f(w)\\
f(rv)=\bar{r}f(v)
\end{array}&r\in\bbbc,\;v,w\in V.
\end{array}$$

 If $A$ is an
algebra over $\bbbc$, $f:A\longrightarrow A$ is called {\it semi-linear involution} if $f$ is
semi-linear  and satisfies (\ref{in}).

\begin{exa}\label{example}
{\em Let $\nu\geq1$ and  ${\bf q}=(q_{i,j})_{1\leq i,j\leq \nu}$ be a $\nu\times\nu-$matrix such
that
\begin{equation}\label{quantum}
q_{i,i}=1\hbox{\;\;and\;\;}q_{i,j}=q_{j,i}=\pm1,1\leq i\neq
j\leq\nu.
\end{equation}
Let $\aa$ be  the unital associative  algebra over $\bbbc $ defined by  generators
$\{t_i^{\pm1}\}_{i=1}^\nu$ and relations
\begin{equation}\label{y3}
\begin{array}{c}
t_it_j=q_{i,j}t_jt_i,\;\;\;1\leq i\neq j\leq\nu.
\end{array}
\end{equation} The algebra $\aa$ is a {\it quantum torus} [BGK, $\S2$].
We have
\begin{equation}\begin{array}{c}
\displaystyle{\aa=\bigoplus_{\sg\in\zn}{\bbbc}t^\sg=\bigoplus_{\sg\in\zn}
({\bbbr}t^\sg\oplus{\bbbr}it^\sg)}\;\hbox{ where}\vspace{1mm}\\
 t^\sg=t_{1}^{n_1}\ldots
t_{\nu}^{n_\nu}\hbox{\; for\; } \sg=(n_1,\ldots,n_\nu)\in\zn.
\end{array}\end{equation}
  There is a semi-linear involution $\ ^-$
on $\aa$ defined by
$$\bar{t_i}=t_i,\quad 1\leq i\leq\nu  \andd \overline{xt^\sg}=
\bar{x}\bar{t^\sg};\quad x\in\bbbc,\;\;\sg\in\zn.$$
For $\sg=(n_1,\ldots,n_\nu)\in\zn$ put $\displaystyle{\kappa_\sg=\prod_{1\leq
i<j\leq\nu}(q_{i,j})^{n_in_j}}$. Let us call $\sg\in\zn$ {\it even} if $\kappa_\sg$ is positive and
call it {\it odd} if $\kappa_\sg$ is negative.  It is easy to see that
\begin{equation}
\kappa_\sg=\kappa_{-\sg},\;t^\sg t^{-\sg}=\kappa_\sg,\;\overline{t^\sg}=\kappa_\sg
t^\sg;\;\;\;\;\sg\in\zn. \label{k3}
\end{equation} Now  consider $\aa$ as an $\bbbr-$algebra and let
$\ell\in\bbbz_{\geq2}$. Put $$E=\left(\begin{array}{cc}
0 & I_\ell \\
-I_\ell & 0
\end{array}\right).$$  $E$ is an invertible
$2\ell\times2\ell$ matrix  and $E^{-1}=-E$. Let  $^*$ be the $\bbbr-$linear automorphism  of
$M_{2\ell}(\aa)$ (the set of $2\ell\times2\ell$ matrices over $\aa$) defined by
$$^*:M_{2\ell}(\aa)\longrightarrow M_{2\ell}(\aa);\;\;X\mapsto E^{-1}\bar{X}^tE,$$ where for $X\in
M_{2\ell}(\aa),$ $X^t$ denotes the transposition of $X$. For  $A=(a_{i,j})_{i,j},
B=(b_{i,j})_{i,j}\in M_{2\ell}(\aa)$, we have
\begin{equation}
(AB)^*=B^*A^*\andd (A^*)^*=A. \label{e1}
\end{equation}
That is,   $^*$ is an involution of $M_{2\ell}(\aa).$  Set $$\bb =\{X\in M_{2\ell}(\aa)\mid
X^*=-X\}.$$ It is easy to check that  $X\in\bb$ if and only if
$$X=\left(\begin{array}{cc}
  A & S  \\
  T & -\bar{A}^t
  \end{array}\right)\hbox{ where } \;A,S,T\in M_{\ell}(\aa)\andd
 \bar{S}^t=S,\;\bar{T}^t=T.$$
 By (\ref{e1}),  $\bb$ is a subalgebra of
 $\mathfrak{gl_{2\ell}}(\aa).$

For  $1\leq r,s\leq2\ell,$ take
 $e_{r,s}$  to be the matrix having $1$ in $(r,s)$ position and
$0$ elsewhere.  Next set $\dot{\hh}:=\sum^{\ell}_{ r=1}\mathbb{R}\dot{h}_r$ where
$\dot{h}_r=e_{r,r}-e_{\ell+r,\ell+r}$ for $1\leq r\leq\ell$.  For $1\leq r\leq\ell$ define
$\ve_r\in\hd^\star$ by $\ve_r(\dot{h}_s)=\d_{r,s}$, $1\leq s\leq \ell$. We can see
\begin{equation}\label{kn1}
\bb=\sum_{\dot{\a}\in\hd^\star}\bb_{\dot{\a}}\;\hbox{where   }
 \;\bb_{\dot{\a}}=\{x\in\bb\mid
[h,x]=\dot{\a}(h)x,h\in\dot{\hh}\},\;\adot\in\dot{\hh}^\star.
\end{equation}
Set $\dot{R}:=\{\dot{\a}\in\hd^\star\mid \bb_{\dot{\a}}\not=\{0\}\}$. Then
$\dot{R}\setminus\{0\}=\rds\cup\rdl$ where
\begin{equation}
\begin{array}{c}
\rds=\{\pm(\ve_r\pm\ve_s)\mid 1\leq r\not=s\leq\ell\}\andd\rdl=
\{\pm2\ve_r\mid 1\leq r\leq\ell\} .\\
\end{array}\label{kn2}
\end{equation}Moreover  we
have
\begin{equation}\label{bb}
\begin{array}{c}
\bb_{\ve_r-\ve_s}=\{ae_{r,s}-\bar{a}
e_{\ell+s,\ell+r}\mid a\in\mathcal{A}\},\;\;r\neq s\vspace{1mm}\\
\bb_{\ve_r+\ve_s}=\{ae_{r,\ell+s}+\bar{a}
e_{s,\ell+r}\mid a\in\mathcal{A}\},\vspace{1mm}\\
\bb_{-\ve_r-\ve_s}=\{ae_{\ell+r,s}+\bar{a}
e_{\ell+s,r}\mid a\in\mathcal{A}\},
\end{array}
\end{equation}
and \begin{equation} \bb_0=\left\{\left(\begin{array}{cc}
  A& 0 \\
  0 & -\bar{A}^t \\
 \end{array}\right)\mid A\hbox{ is diagonal }
\right\}.\label{b0}
\end{equation}

Now we want to define a real  form $\fm$ on $\bb$. For this, we
first define $\epsilon :\aa\longrightarrow \mathbb{R}$ by linear
extension of
$$\epsilon(t^\sg)=\left\{\begin{array}{cc}
  1 & \hbox{if\ }\sg=0 \\
  0  & \hbox{if\ }\sg\neq 0
  \end{array}\right. \andd\epsilon(it^\sg)=0.$$
Then $(a,b)\mapsto \epsilon(ab)$ is a  non-degenerate symmetric real bilinear form on $\aa$
preserved by   $\ ^-$.  Now one can define a symmetric invariant real bilinear form on $M_{2\ell}(\aa)$ (and so on
$\mathfrak{gl_{2\ell}}(\aa)$) as follows:

$$(A, B)=\epsilon(\tr(AB));\;\;\;A,B\in
M_{2\ell}(\aa).$$
For $1\leq r,s\leq \ell$, we have
$(\dot{h}_r,\dot{h}_s)=\epsilon(\tr(\dot{h}_r\dot{h}_s))=\ep(2\d_{r,s})= 2\ve_{r}(\dot{h}_s)$, so
the form restricted to $\dot{\hh}$ is non-degenerate. Also
\begin{equation}
t_{\ve_r}=\frac{1}{2}\dot{h}_r\andd (t_{\ve_r},t_{\ve_s})=\frac{1}{2}\d_{r,s};\;\;1\leq r,s\leq
\ell. \label{e6}
\end{equation} It is easy to check that
\begin{equation}
\hbox {the form restricted to $\bb$ is non-degenerate, symmetric
and  invariant}. \label{e'5}
\end{equation}

Next, we would like to make $\bb$ to a $\zn-$graded Lie algebra. As in [AABGP, \S III.3], we start
with a gradation on $M_{2\ell}(\aa)$ (as a vector space over $\bbbr$). For $1\leq p,q\leq 2\ell$, set
$$\begin{array}{c}\hbox{deg}(it^\sg
e_{pq})=\hbox{deg}(t^\sg e_{pq})=\sg.
\end{array}$$
This defines a $\zn-$grading on $M_{2\ell}(\aa)$ and in turn on
$\mathfrak{gl_{2\ell}}(\aa)$. Moreover, using (\ref{k3}), we have
$(t^\sg e_{p,q})^*=\kappa_\sg E^{-1}t^\sg e_{q,p}E$ for $1\leq
p,q\leq2\ell$ and $\sg\in\zn$. Thus   the involution $^*$
preserves the grading on $M_{2\ell}(\aa)$ and so
\begin{equation}\label{kn3}
\bb\;\hbox{is  a $\zn-$graded subalgebra of}\; \mathfrak{gl_{2\ell}}(\aa).
\end{equation}
 It is easy to check
that the form on $\bb$ is a graded form. Also it
 can be easily seen from (\ref{bb}) and (\ref{k3}) that
\begin{equation}
\bb_{\dot{\a}}=\bigop_{\sg\in\bbbz^\nu }(\bb_{\dot{\a}}\cap\bb^\sg);\;\; \dot{\a}\in \dot{R}.
\label{e3}
\end{equation}

Next  let $\gg$  be the derived algebra of $\bb.$ Then
\begin{equation*}
\gg=\gg_0+\sum_{\adot\in\dot{R}^\times}\gg_{\adot} \label{b01}
\end{equation*}
where  $\gg_{\adot}=\bb_{\adot}$ for $\adot\in\dot{R}^\times$ and
$\gg_0=\bb_0\cap\gg.$ Now let $1\leq r\neq s\leq\ell$ and
$\sg\in\zn$. Using (\ref{k3}) and (\ref{bb}), we have
\begin{equation}\label{sbb}
\begin{array}{c}
\gg_{\ve_r-\ve_s}^\sg=\bbbr t^\sg(e_{r,s}-\kappa_\sg
e_{\ell+s,\ell+r})
+i\bbbr t^\sg(e_{r,s}+\kappa_\sg  e_{\ell+s,\ell+r}),\vspace{1mm}\\
\gg_{\ve_r+\ve_s}^\sg=\bbbr t^\sg(e_{r,\ell+s}+\kappa_\sg
e_{s,\ell+r}) +i\bbbr t^\sg(e_{r,\ell+s}-\kappa_\sg e_{s,\ell+r})
,\vspace{1mm}\\
\gg_{-\ve_r-\ve_s}^\sg=\bbbr t^\sg(e_{\ell+r,s}+\kappa_\sg
e_{\ell+s,r}) +i\bbbr t^\sg(e_{\ell+r,s}-\kappa_\sg e_{\ell+s
,r}) ,\vspace{1mm}\\
\gg_{2\ve_r}^\sg=\left\{\begin{array}{cc}
\bbbr t^\sg e_{r,\ell+r}& \hbox{if $\sg$ is even}\\
i\bbbr t^\sg e_{r,\ell+r}& \hbox{if $\sg$ is odd},
\end{array}\right.
\vspace{1mm}\\
\gg_{-2\ve_r}^\sg=\left\{\begin{array}{cc}
\bbbr t^\sg e_{\ell+r,r}& \hbox{if $\sg$ is even}\\
i\bbbr t^\sg e_{\ell+r,r}& \hbox{if $\sg$ is odd}.
\end{array}
\right.
\end{array}
\end{equation}

For $\sg=(n_1,\ldots,n_\nu), \tau=(m_1,\ldots,m_\nu)\in\zn,$
define
\begin{equation*}
\begin{array}{c}
g_\sg^\tau=\displaystyle{\prod_{1\leq i\leq j\leq
\nu}q_{i,j}^{n_im_j}}\andd f_\sg^\tau:=g_\sg^\tau g_\tau^\sg.
\end{array}
\end{equation*}
One can easily see that  for $\sg,\tau,\gamma\in\zn,$
\begin{equation}
\begin{array}{c}\kappa_\sg=g_\sg^\sg,\;
f_\sg^\tau\kappa_\sg\kappa_\tau=\kappa_{\sg+\tau},\;
g_{\sg+\tau}^\gamma=g_\sg^\gamma g_\tau^\gamma,\;
g_\gamma^{\sg+\tau}=g_\gamma^\sg g_\gamma^\tau. \label{rev1}
\end{array}
\end{equation}
 Now for
$a,b\in\bbbr,$ $1\leq r\neq s\leq\ell$ and $\sg,\tau\in\zn,$ set

$$\begin{array}{l}
A_{a,b}^{\sg,r,s}:=at^\sg(e_{r,s}-\kappa_\sg
e_{\ell+s,\ell+r})+ibt^\sg(e_{r,s}+\kappa_\sg e_{\ell+s,\ell+r}),\\
B_{a,b}^{\sg,r,s}:=at^\sg(e_{r,\ell+s}+\kappa_\sg
e_{s,\ell+r})+ibt^\sg(e_{r,\ell+s}-\kappa_\sg e_{s,\ell+r}),\\
C_{a,b}^{\sg,r,s}:=at^\sg(e_{\ell+r,s}+\kappa_\sg
e_{\ell+s,r})+ibt^\sg(e_{\ell+r,s}-\kappa_\sg
e_{\ell+s,r}),\\
m_{r}^{\sg,\tau}:=e_{r,r}-\kappa_{\sg+\tau}e_{\ell+r,\ell+r},\hspace{.8cm}
n_{r}^{\sg,\tau}:=e_{r,r}+\kappa_{\sg+\tau}e_{\ell+r,\ell+r}.\end{array}$$

Using (\ref{rev1}), we have

\begin{eqnarray}
\;[A_{a,b}^{\sg,r,s},A_{c,d}^{\tau,s,r}]&=&
g_{\sg}^{\tau}((ac-bd)(m_{r}^{\sg,\tau}-f_{\sg}^{\tau}
m_{s}^{\sg,\tau})\nonumber\\
&+&i(ad+bc)(n_{r}^{\sg,\tau}-f_{\sg}^{\tau}n_{s}^{\sg,\tau}))t^{\sg+\tau},\label{rev2}\\
\;[B_{a,b}^{\sg,r,s},C_{c,d}^{\tau,r,s}]&=&g_\sg^\tau((ac+bd)\kappa_\tau(m_{r}^{\sg,\tau}
+f_\sg^\tau\kappa_{\sg+\tau}m_{s}^{\sg,\tau})\nonumber\\
&+&i(bc-ad)\kappa_\tau(n_{r}^{\sg,\tau}-f_\sg^\tau\kappa_{\sg+\tau}
n_{s}^{\sg,\tau}))t^{\sg+\tau}.\nonumber
\end{eqnarray}
Also

 \begin{equation} [t^\sg e_{r,\ell+r}, t^\tau
e_{\ell+r,\ell+r}]=\left\{\begin{array}{ll}
g_{\sg}^\tau m_{r}^{\sg,\tau}t^{\sg+\tau}& \hbox{if}\;\;\kappa_\sg\kappa_\tau=1\\
g_{\sg}^\tau n_{r}^{\sg,\tau}t^{\sg+\tau}&
\hbox{if}\;\;\kappa_\sg\kappa_\tau=-1.
\end{array}\right.\label{rev3}
\end{equation}
We recall   that  for $1\leq r\leq \ell,$
$\dot{h}_{r}=e_{r,r}-e_{\ell+r,\ell+r}$ and put
$\ddot{h}_r:=e_{r,r}+e_{\ell+r,\ell+r}.$ Now let $\gamma\in\zn.$
We say $\gamma$ satisfies {\it even property} (resp. {\it odd
property}) if there exists $\sg,\tau\in\zn$ such that
$\sg+\tau=\gamma$ and $\kappa_\sg\kappa_\tau=1$ (resp.
$\kappa_\sg\kappa_\tau=-1 $).  Using (\ref{rev2}) and
(\ref{rev3}), we have
\begin{eqnarray*}
\gg_0^\gamma&=&\displaystyle{ \sum_{\adot\in\rd^\times}
\sum_{\hbox{\tiny$\begin{array}{c}\tau,\sg\in\zn\\
\sg+\tau=\gamma
\end{array}$}}}[\gg_{\adot}^\sg,\gg_{-\adot}^\tau]=\sum_{1\leq r\neq s\leq\ell}\gg_0^{\gamma,r,s}\\
\end{eqnarray*}
where for $1\leq r\neq s\leq\ell,$ $\gg_0^{\gamma,r,s}$ equals
$$\left\{\begin{array}{ll}
(\bbbr\dot{h}_r+i\bbbr\ddot{h}_r)t^\gamma&\hbox{if $\gamma$ is
even and satisfies odd property,}\\
(\bbbr\dot{h}_r+i\bbbr(\ddot{h}_r-\ddot{h}_s))t^\gamma&\hbox{if
$\gamma$ is
even and dose not satisfy odd property,}\\
(\bbbr\ddot{h}_r+i\bbbr\dot{h}_r)t^\gamma&\hbox{if $\gamma$ is
odd and satisfies even property,}\\
(\bbbr(\ddot{h}_r-\ddot{h}_s)+i\bbbr\dot{h}_r)t^\gamma&\hbox{if
$\gamma$ is
odd and dose not satisfy even property.}\\
\end{array}
\right.$$

Now by mimicking the argument given in [AABGP, \S III.4] we can
check that $(\gg,(\cdot,\cdot)|_{\gg\times\gg},\dot\hh)\in\dd$.
Therefore using Theorem \ref{cnonstruct}, we can add to $\gg$ some
central elements and some derivations to obtain a Lie algebra
$\LL\in \T$ which is tame and  of type $C_\ell$. However $\dot\hh$
is a proper subset of $\gg_0\cap\gg^0$  and so $\LL$ does not
satisfy the second axiom of an EALA (see (\ref{rs})). Also by
(\ref{rs}) and (\ref{sbb}),
 $\dim(\LL_{\pm\ve_r\pm\ve_s+\sg})=\dim(\gg_{\pm\ve_r\pm\ve_s}^\sg)=2$ for
 $1\leq r\neq s\leq\ell$ and $\sg\in\zn.$ This phenomena never  happens
 for an EALA [AABGP, Theorem I.1.29].}
\end{exa}
\section{ Nullity zero}
It is known that a Lie algebra  $\LL$ over $\bbbc$ is a finite
dimensional simple Lie algebra if and only if
 $\LL$ is a finite dimensional tame EALA if and only if it is tame EALA of nullity 0. In this section, we first  study  the elements  in $\T$   of nullity
zero. We find a necessary and sufficient condition  for   a Lie
algebra over $\bbbf$ to be an element in $ \T$  of nullity zero.
Then we   find a new characterization for
 finite dimensional simple Lie algebras  containing  a nonzero
 split
 torus. In fact we prove that a
 finite dimensional simple Lie algebra  containing  a nonzero
 torus is an element of $\T$ satisfying certain properties.
Conversely, an element of $\T$  satisfying these properties is a
 finite dimensional simple Lie algebra  containing  a nonzero
 torus.

\begin{DEF}
{\em Let $\D$  be an irreducible finite root system
 and $G$ be an abelian group. Let $\LL=\op_{\mu\in\D}\op_{g\in G}\LL^g_\mu$ be a $(\D,G)-$graded Lie
algebra with grading pair $(\fg,\fh)$. We call $\LL$  {\it
strictly   division $(\D,G)-$graded Lie algebra} if for each
$\mu\in\D^\times$, $g\in G$ and  $0\neq x\in\LL^g_\mu$, there
exists $y\in\LL^{-g}_{-\mu}$ such that $[x,y]=t_\mu$.}
\end{DEF}
\begin{thm}
Let  $\D$ be an irreducible finite root system and
$(\aa,[.,.]_\aa)$ be a strictly division  $\D-$graded Lie algebra
with grading pair $(\fg,\fh)$. Let $(E, [.,.]_E)$ be a Lie algebra
for which there exist linear maps
\begin{equation}
\begin{array}{l}
\rho:E\longrightarrow (\hbox{\em Der}(\aa))^\fg=\{d\in\hbox{\em Der}(\aa)\mid d(\fg)=0\}\\
\tau:E\times E\longrightarrow \aa \end{array} \label{maps} \end{equation} satisfying the following
conditions:
\begin{equation}
\begin{array}{l}
\tau(x,y)=-\tau(y,x),\\
\hbox{\em ad}(\tau(x,y))=[\rho(x),\rho(y)]-\rho([x,y]_E),\\
\tau([x,y]_E,z)+\tau([z,x]_E,y)+\tau([y,z]_E,x)=\\
\rho(x)(\tau(y,z))+\rho(z)(\tau(x,y))+\rho(y)(\tau(z,x)),
\end{array}
\label{con}
\end{equation}
for $x,y,z\in E$.  Set $\LL:=\aa\op E$ and define an
anti-commutative product $[.,.]$ on $\LL$ as follows:
\begin{equation}
\begin{array}{l}
\ [a,b]=[a,b]_\aa,\\
\ [x,a]=-[a,x]=\rho(x)(a),\\
\ [x,y]=[x,y]_E+\tau(x,y),\end{array} \;\; a,b\in\aa,\;x,y\in E.\label{mult}\end{equation}  Then
$(\LL,[.,.])$ is a Lie algebra containing $\aa$ as an ideal such that $\LL/\aa\cong E.$
 Moreover if $W$ is a finite dimensional subspace  of $Z(\LL)$ and $\LL$ is equipped with a
 non-degenerate   symmetric invariant bilinear form $(\cdot,\cdot)$ which satisfies the following conditions:
\begin{equation}
\begin{array}{l}
(\cdot,\cdot)|_{\fg\times\fg} \;\hbox{ is nontrivial,}\;
 (\fh,W)=0,\\
 (\cdot,\cdot)|_{W\times W} \;\hbox {is non-degenerate,}\\
\end{array}\label{p1c}
\end{equation}
   then
 $(\LL, (\cdot,\cdot),\fh\op W)$ is an element   of $\T$ of nullity zero. Conversely if $\LL\in\T$ is
 of  nullity zero, then there exist a strictly division $\D-$graded ideal $\aa$  of $\LL$ with grading pair $(\fg,\fh)$ and a  vector
subspace  $E\subseteq C_\LL(\fg)$ such that $\LL=\aa\op E$. In
this case $E$ is a Lie algebra isomorphic to $\LL/\aa,$ and there
exist linear maps $\rho$ and $\tau$ as in (\ref{maps}) satisfying
(\ref{con}) and (\ref{mult}).
  \label{nz}
\end{thm}
\noindent{\bf Proof.} One can easily check that $\LL=\aa\op E$ is
a Lie algebra. Since $[\fh,E]\subseteq[\fg,E]=\{0\},$ we have a
weight space decomposition of $\LL$ as $\LL=\op_{\a\in\D}\LL_\a$
where $\LL_0=\aa_0\op E$ and $\LL_\a=\aa_\a$ for $\a\in
\D^\times$. Suppose now that $(\cdot,\cdot)$  is a form on $\LL$
and  $W$ is a subspace of $Z(\LL)$ satisfying (\ref{p1c}). We
extend any $\a\in \hh^\star$ to a linear functional of $\hh\op W$
by defining $\a|_W=0.$ Since the form on $\fg$ is nontrivial,
$(\cdot,\cdot)|_{\fg\times\fg}$ is a nonzero scalar multiple of
the Killing form (see the proof of  Lemma \ref{ac}) and so
$(\cdot,\cdot)|_{\fh\times\fh}$ is non-degenerate, therefore the
form on $\fh\op W$ is non-degenerate. Now since $W\subseteq
Z(\LL),$ we have $\LL_\a=\{x\in\LL\mid [h,x]=\a(h)x,h\in\fh\op
W\},$ $\a\in \D$. Thus  $(T1)$ and $(T2)$ hold. Also we have
$(T3)$ since $\LL$ is a strictly division $\D-$graded Lie algebra
and the form on $\fg$ is  a nonzero scalar multiple of the Killing
form such that $(\fh,W)=0$. Finally, we have $(T4)$, $(T5)$ and
$(T6)$ because $\D$ is an irreducible finite root system.

Conversely, suppose $\LL$  is an element of $\T$ of nullity zero.
Then $\LL=\op_{\a\in\rd}\LL_{\a}$ where $\rd$ is an irreducible
finite root system. Using $(T3)$ together with Theorem
\ref{deltagraded}, one can get that $\aa:=\LLc$ is a strictly
division  $\rd-$graded Lie algebra with grading pair
$(\dg,\dot\hh)$ where $\dgg$ and $\dot\hh$ are defined as in
Theorem \ref{m12}. Consider a subspace $D$ of $\LL$ such that
$\LL=D\op\aa$ and let $x\in D$. Since $\LLc$ is an ideal, $\ad
(x)|_{\LLc}\in \hbox{Der}(\LLc)$. Using the complete reducibility
of $\LLc$ as a $\dot\gg-$module and the first Whitehead lemma for
$\dgg-$modules, we can apply [Be, Proposition 3.2] to each element
of a basis of $D$ and in this way construct  a subspace $E$ such
that $\LL=\LLc\op E$ and $[E,\dgg]=0.$ So we have a linear map
$\rho:E\longrightarrow \hbox{Der}(\aa)^{\dgg}$ given by
$\rho(x)=\ad(x)|_{\LLc}$ for $x\in E$. Since $\LL$ is a Lie
algebra, the Jacobi identity together with the anti-commutative
product on $\LL$ implies  the existence of a Lie bracket on $E$
and a linear map $\tau$ as in (\ref{maps}) satisfying (\ref{con})
and (\ref{mult}).\qed
\smallskip

It is a well-known fact that tame EALAs of nullity zero are
 finite dimensional, but   a tame element in
 $\T$ of nullity zero need not to be finite dimensional. More precisely suppose that $\gg$ is a finite dimensional
 split simple Lie
algebra  over a field $\bbbf$ of characteristic zero with root
system $\D$, splitting  Cartan subalgebra $\hh$ and the Killing
form $\kappa(\cdot,\cdot).$ Let $\bbbk$ be a field extension of
$\bbbf$. Since $\gg$ is central simple,
$\bar{\gg}:=\gg\ot_{\bbbf}\bbbk$ with Lie bracket $[x\ot a, y\ot
b]=[x,y]\ot ab,$ $x,y\in \gg,\; a,b\in \bbbk$ is a simple Lie
algebra over $\bbbk$ and consequently $\bar{\gg}$ is centerless.
 For $\a\in\D^\times$ and  the corresponding
root space  $\gg_\a,$ fix $x_{\pm\a}\in\gg_{\pm\a},$ such that
$[x_\a,x_{-\a}]=h_\a,$ where $h_\a=2t_\a/\kappa(\a,\a).$  One can
easily see that $\bar{\gg}$ has a weight space decomposition  with
respect to $\hh\ot 1$ as $\bar{\gg}=\op_{\a\in\D}\bar{\gg}_\a,$
where $\bar{\gg}_\a:=\gg_\a\ot_\bbbf \bbbk.$ Also for each
$a\in\bbbk\setminus \{0\},$ there exists $b\in \bbbk$ such that
$[x_\a\ot a , x_{-\a}\ot b]=h_\a\ot 1.$ Therefore $\bar{\gg}$ is a
strictly division $\D-$graded Lie algebra with grading pair
$(\gg\ot1,\hh\ot1).$ Now let $\bbbk$ be equipped with a
non-degenerate symmetric bilinear form $f(\cdot,\cdot):\bbbk\times
\bbbk\longrightarrow \bbbf$ which is invariant in the sense that
\begin{equation}
f(ab,c)=f(a,bc);\;\; a,b,c\in \bbbk. \label{rev}
\end{equation}
Then
\begin{equation}
(\cdot,\cdot):\bar{\gg}\times\bar{\gg}\longrightarrow\bbbf;\;\;
(x\ot a, y\ot b)\mapsto \kappa(x,y)f(a,b),
\end{equation}
is a non-degenerate symmetric invariant bilinear form  on
$\bar{\gg}$ as a Lie algebra over $\bbbf.$ Therefore by  Theorem
\ref{con}, $\bar{\gg}$ is a Lie algebra belonging to $\T$ of
nullity zero which is tame as $\LL_c=\LL$. Now if $[\bbbk:\bbbf]$
is infinite, then $\bar{\gg}$ is an infinite dimensional Lie
algebra belonging to $\T$ which is tame and of nullity zero.
\begin{exa}{\rm
Let $\bbbf$ be the field of rational numbers $\bbbq$ and $\aa$ be
the set of all positive square-free integers. Set
$\mathcal{B}:=\{\sqrt a | a\in \aa\}.$  Then $\mathcal{B}$ is a
linearly independent set over  $\bbbq$ (see e.g.  [W]). Therefore
$\bbbk:=\hbox{span}_\bbbq \bb$ is a field extension of $\bbbq$ and
$\bb$ is a basis of $\bbbk$ over $\bbbq.$ Now define $f(\sqrt
a,\sqrt b):=\d_{a,b}a,\; a,b\in\bbbk,$ where $\delta_{a,b}$
denotes the Kronecker delta. Extend this linearly to a bilinear
form $f(\cdot,\cdot):\bbbk\times \bbbk\longrightarrow \bbbq$. It
is easy to  see that $f$ is non-degenerate,  symmetric and
invariant and so by the above discussion $\bar{\gg}$ is an
infinite dimensional Lie algebra belonging to $\T$ which is tame
and of nullity zero.}
\end{exa}
We recall that for $(\LL,(\cdot,\cdot),\hh)\in\T$, we have
$\LL_0=\hh\op\hh^\bot$, where  $\hh^\bot$ is the orthogonal
complement of $\hh$ in $\LL_0.$
\begin{lem}
Let $\LL\in\mathcal{T}$ satisfy
$\sum_{\a\in\rtimes}\hh_\a=\hh^\bot$ (see (\ref{new})). If the
nullity of $\LL$ is zero, then $\LLc$ is a simple Lie algebra and
$\LL=\LLc\odot Z(\LL)$ where $\odot$ means orthogonal direct sum.
Moreover  if $\LL$ is tame, then  $\LL=\LL_c$ is a simple Lie
algebra. \label{gs}
\end{lem}
\noindent{\bf Proof.}  We know from  (\ref{G2}) and Lemma \ref{m3}
that $\LL_0=\hh\op\hh^\bot=\dot\hh\op\hh^0\op D\op
W\op\sum_{\a\in\rtimes}\hh_\a$ and $W\subseteq Z(\LL)$. Since the
nullity of $\LL$ is  zero, $\hh^0\op D=\{0\}$. Also by Proposition
\ref{m1}($iv$), we have $\sum_{\a\in\rtimes}\hh_\a\subseteq\LLc,$
so by Proposition \ref{m1}($viii$), $\LL=\LLc\odot W$. We know
that the forms on $\LL$ and on $W$ are non-degenerate, so the form
on $\LLc$
 is non-degenerate. Therefore by Proposition \ref{m1}($viii$), $\LLc$  is
centerless. Hence $W=Z(\LL)$ and by Corollary \ref{c02}, $\LLc$ is
a simple $\rd-$graded  Lie algebra. If moreover  $\LL$ is tame,
then by Lemma \ref{m3}($iii$) we have $W=0$ and so $\LL=\LLc$  is
a simple $\rd-$graded Lie algebra. \qed

\begin{thm}
Let $\LL$ be a finite dimensional Lie algebra. Then $\LL$ is a
tame  element in $ \T$  with $\sum_{\a\in\rtimes}\hh_\a=\hh^\bot$
if and only if $\LL$ is simple Lie algebra containing a nonzero
 torus.
\end{thm}
\noindent{\bf Proof.} Let $\LL\in \mathcal{T}$ be finite
dimensional and tame which satisfies $\sum_{\a\in\rtimes}\hh_\a$
$=\hh^\bot$. By [AKY, (2.1)], $R^0$ is a semilattice, so the
nullity of $\LL$ is zero since $\LL$ is finite dimensional. Now
Lemma \ref{gs} implies that $\LL$ is a finite dimensional simple
Lie algebra. Moreover $\LL$ contains  a nonzero  torus since $\LL$
satisfies $(T2)$.

Conversely, let $\LL$ be a finite dimensional simple Lie algebra
containing  a nonzero  torus. Since $\bbbf$ is of characteristic
zero, the Killing form on $\LL$ is non-degenerate, so $\LL $
satisfies $(T1)$. Let $0\neq \hh$ be a maximal torus, then by [Se,
Corollary to Lemma I.1.6], the form on $\hh$ is non-degenerate, so
$ (T2)$ holds. By [Se, \S I.1 and \S III.1], $\LL=\op_{\a\in
\rd}\LL_{\adot}$, where $\rd$ is an irreducible finite root
system. Now we have $(T3)$ by [Se, Lemma I.1.3 and Corollary to
Lemma I.1.4]. The axioms  $(T4)$, $(T5)$ and $(T6)$ are satisfied
because $\dot R$ is an irreducible
 finite root system. Finally, since $\LL$ is  simple, we have
$\LL_0=\sum_{\a\in\rd^\times}[\LL_\a,\LL_{-\a}]$ then  $\LL$ is
tame and $\sum_{\a\in\rtimes}\hh_\a=\hh^\bot$.
 \qed

\begin{exa}
{\rm Let $\bbbf=\bbbr$. Each noncompact real Lie algebra  is a
tame finite dimensional  element of $\mathcal{T}$ with toral
subalgebra $\hh$ and root system $R$  satisfying
$\sum_{\a\in\rtimes}\hh_\a=\hh^\bot$(see [Se, \S V.8]).}
\end{exa}
\section{ Appendix}
In [AKY], the authors axiomatically     introduced a class of Lie
algebras called {\it toral type extended affine Lie algebras}. The
root system of a toral type extended affine Lie algebra is an EARS
(see Definition \ref{ears}) [AKY, \S1]. One can see that the
axioms of a toral type extended affine Lie algebra are slightly
different from $(T1) - (T6)$. In this section we show that a Lie
algebra $\LL\in \T$ is a toral type extended affine Lie algebra.
For this, it is enough to prove that for a Lie algebra $\LL\in\T$
with root system $R$, the $\bbbq-$span of the root system $R,$
$\v_\bbbq$, is a finite dimensional vector space and $R$ is a
discrete subset of $\bbbr\ot_\bbbq \v_\bbbq.$
\begin{DEF}\label{ears}
{\em Let $V$ be a nontrivial finite dimensional real vector space
with a positive semidefinite symmetric bilinear form
$(\cdot,\cdot)$. Let $R$ be a  nonempty subset of $V.$ Let
$$
R^{\times}=\{\a\in R\mid (\a,\a)\not=0\}\quad\hbox{and} \quad
R^{0}=\{\a\in R\mid (\a,\a)=0\}.
$$
Then $R=R^{\times}\uplus R^{0}$.  $(R,(\cdot,\cdot))$ (or $R$ if
there is  no any confusing) is called   an {\it extended affine
root system }  in $V$ if  $R$ satisfies the following 5 axioms:

\noindent (R1) $R=-R$.

\noindent (R2)$\;R \hbox{ spans }V.$

\noindent (R3)$\;$ $R$ is discrete in $V.$

\noindent (R4)$\;$ For  $\a\in R^{\times}$ and $\b\in R$, there
exist $d, u\in {\mathbb Z}_{\geq 0}$ such that
$$\{ \b\hbox{\tiny{+}} n\a\mid n \in {\bbbz} \} \cap R =
\{ \b-d\a,\ldots,\b\hbox{\tiny{+}} u\a \} \hbox{ and
}d-u=2{(\b,\a)}/{(\a,\a)}.
$$
\noindent (R5)(a) $R^{\times}$ cannot be written as a disjoint
union of two nonempty subsets which are orthogonal with respect to
the form.

(b) For any $\d\in R^{0}$, there exists $\a\in R^{\times}$ such
that $\a+ \d\in R$.

$R$ is called {\it reduced} if it satisfies:

\noindent (R6)$\;\a\in R^{\times} \Rightarrow 2\a\not\in R$.}
\end{DEF}

Now let $\LL\in\T$ with root system $R$. Let
$\v_\bbbf:=\hbox{span}_\bbbf(R)$ and $(\cdot,\cdot)_\bbbf$ be the
form
 restricted to $\v_\bbbf$. Take  $\v_\bbbf^0$ to be the radical of
the
 form $(\cdot,\cdot)_\bbbf$ and
 $\overline{\v}_\bbbf={\v_\bbbf}/{\v_\bbbf^0}$.
 Let $^-:\v_\bbbf\longrightarrow\overline{\v}_\bbbf$ be the
 canonical map. The  form $(\cdot,\cdot)_\bbbf$ on $\v_\bbbf$
 induces a unique non-degenerate symmetric bilinear form $(\cdot,\cdot)_\bbbf$
on
 $\overline{\v}_\bbbf$ so that
 $(\bar{\a},\bar{\b})_\bbbf=(\a,\b)_\bbbf$, $\a,\b\in\v_\bbbf.$
Let $\bar{R}$  be the image of $R$ under $^-.$ For  $\a\in R^0,$
Proposition \ref{m1}($i$) gives that $(\a, R)=0$ so
 $(\a,\v_\bbbf)_\bbbf=0.$  It means that  $\a\in\v_\bbbf^0.$ Also if
$\a\in\v_\bbbf^0\cap
 R$, then $(\a,\a)=0$, in other words $\a\in R^0$.  Therefore we
 have
 \begin{equation}\label{q1}
 R^0=R\cap\v_\bbbf^0\andd \bar{R}\setminus\{0\}=\{\bar{\a}\mid\a\in \rtimes\}.
 \end{equation}
Let $\v_\bbbq=\hbox{span}_\bbbq(R)$ and  $(\cdot,\cdot)_\bbbq$ be
the
 form
 restricted to $\v_\bbbq$. Let  $\v_\bbbq^0$  be the radical of the
 form $(\cdot,\cdot)_\bbbq$ and
 $\tilde{\v}_\bbbq={\v_\bbbq}/{\v_\bbbq^0}$.
 Let $\tilde{}:\v_\bbbq\longrightarrow\tilde{\v}_\bbbq$ be the
 canonical map. The  form $(\cdot,\cdot)_\bbbq$ on $\v_\bbbq$
 induces a unique non-degenerate symmetric bilinear form $(\cdot,\cdot)_\bbbq$
on
 $\tilde{\v}_\bbbq$ so that
 $(\tilde{\a},\tilde{\b})_\bbbq=(\a,\b)_\bbbq$, $\a,\b\in\v_\bbbq.$
Set $\tilde{R}$ to be the image of $R$ under\; $\tilde{}\;.$
 One can easily check that
\begin{equation}
R^0=\v_\bbbq^0\cap R\andd \v_\bbbq^0=\v_\bbbf^0\cap\v_\bbbq.
\label{q2}
\end{equation}
So $\tilde{R}\setminus\{0\}=\{\tilde{\a}\mid \a\in\rtimes\}.$ Now
let $\tilde{\a}\in\tilde{R}\setminus\{0\}$ and define
$w_{\tilde{\a}}\in\GL(\tilde{\v}_\bbbq)$ by:
$$w_{\tilde{\a}}(\tilde{\b})=\tilde{\b}-\frac{2(\tilde{\b},\tilde{\a})_\bbbq}{(\tilde{\a},\tilde{\a})_\bbbq}
\tilde{\a};\;\; \b\in\v_\bbbq.$$ It is easily seen that
\begin{equation}
w_{\tilde{\a}}(\tilde{\b})={(w_\a(\b))}^{\tilde{}}\;;\;\;\a,\b\in
R, \label{q3}
\end{equation}
where ${(w_\a(\b))}^{\tilde{}}$ is the image of $w_\a(\b)$ under
map $\tilde{\;\;}.$
\begin{lem}
$\tilde{\v}_\bbbq$ is a finite dimensional vector space over
$\bbbq$ and $\tilde{R}$ is a  finite root system in
$\tilde{\v}_\bbbq$. \label{aaa}
\end{lem}
\noindent{\bf Proof.} Since $\v_\bbbf$ is a subset of $\hh^\star,$
dim$\v_\bbbf<\infty$ and so $\dim\overline{\v}_\bbbf<\infty.$ Now
using  the same argument as in [AABGP, Lemma I.2.10], one can see
that  $\bar{R}$ is a finite set.
 By (\ref{q2}), we have
$\v_\bbbq^0=\v_\bbbf^0\cap\v_\bbbq.$ So the map
$\phi:\tilde{\v}_\bbbq\longrightarrow\overline{\v}_\bbbf$ defined
by $\phi(\tilde{\a})=\bar{\a}$ is well defined and injective. One
can see that $\phi(\tilde{R})=\bar{R}$ and so $\tilde{R}$ is
finite. Therefore $\tilde{\v}_\bbbq$ is finite dimensional since
$\tilde{\v}_\bbbq=\hbox{span}_\bbbq(\tilde{R})$. Now we are done,
using (\ref{q3}) together with Proposition \ref{v2}.\qed
\\

Now let $\Pi=\{\tilde{\a}_1,\ldots,\tilde{\a}_s\}$  be a base for
$\tilde{R}.$ Fix a preimage $\dot{\a}_i\in R$ of $\tilde{\a}_i$
for $i=1,\ldots,s$ and set
$\dot{\v}_\bbbq:=\hbox{span}_\bbbq\{\dot{\a}_1,\ldots,\dot{\a}_s\}.$
Then $\v_\bbbq=\dot{\v}_\bbbq\op\v_\bbbq^0.$ Define
$\dot{R}:=\{\dot{\a}\in\dot{\v}_\bbbq\mid \adot+\d\in R\;\hbox{for
some} \;\d\in \v_\bbbq^0\} $ and set
$\dot{R}^\times:=\dot{R}\setminus\{0\}.$  $\rd$ is mapped
bijectively to $\tilde{R}$ under $\tilde{\;\;}$. Hence $\rd$ is a
finite root system in $\dot{\v}_\bbbq$ isomorphic to $\tilde{R}.$
For $\adot \in\dot{R}^\times$ let
$S_{\adot}=\{\d\in\v_\bbbq^0\mid\adot+\d\in R\}$, then
\begin{equation}
R=R^0\uplus(\uplus_{\adot\in\dot{R}^\times}(\adot+S_{\adot})),
\label{q4}
\end{equation}
where $\uplus$ means disjoint union.
\begin{lem}
 $\sum_{\adot\in\dot{R}^\times}S_{\adot}\subset R^0$  and $\v_\bbbq^0=\hbox{span}_\bbbq
(R^0)$. \label{aa}
\end{lem}
\noindent{\bf Proof.} Using   the same  argument  as in [AABGP,
Proposition II.2.11(b)], we have $S_{\adot}\subset R^0$ for
$\adot\in\dot{R}^\times$. This then  together with (\ref{q4})
implies the equality in the statement.\qed

\begin{lem}
 $\v_\bbbq=\hbox{span}_\bbbq(R)$ is a
finite dimensional vector space over $\bbbq$.\label{fin}
\end{lem}
\noindent{\bf Proof.}  Axiom  (T6) says that  $\la R^0\ra$ is a
free abelian group of finite rank. If $R^0=\{0\}$, then Lemma
\ref{aa} together with Lemma \ref{aaa} implies that
$\v_\bbbq\simeq\tilde{\v}_\bbbq$ is finite dimensional. Now
suppose that $R^0\neq\{0\}.$ Let $n\in\bbbn\setminus\{0\}$ and
$\{\d_i\}_{i=1}^n$ be a $\bbbz-$basis for $\la R^0\ra.$ Using
Lemma \ref{aa}, we have $\v_\bbbq^0=\hbox{span}_\bbbq\{\d_i\mid
1\leq i\leq n\}$ and consequently  $\dim\v_\bbbq^0<\infty$. Also
by Lemma \ref{aaa},
$\dim\dot{\v}_\bbbq=\dim\tilde{\v}_\bbbq<\infty$. Therefore
$\dim\v_\bbbq=\dim\dot{\v}_\bbbq+\dim\v^0_\bbbq<\infty.$ \qed
\\

Next set $\v:=\bbbr\ot_\bbbq\v_\bbbq.$ Identifying $\v_\bbbq$ as a
subset of $\v$, we can prove the following proposition:
\begin{pro}
$R$ is a discrete subset of $\v.$\label{dis}
\end{pro}
\noindent{\bf Proof.} If $R^0=\{0\}$, then $R=\rd$ and so $R$ is a
discrete subset of $\v$. Now let $R^0\neq\{0\}$ and consider a
$\bbbz-$basis $\{\d_i\}_{i=1}^n$ for $\la R^0\ra$. Using
(\ref{q4})  together with
 Lemma \ref{aa}, we have
 $R\subseteq\op_{j=1}^n\bbbz\d_j\op(\op_{i=1}^s\bbbz\adot_i)\subseteq \v$.
 This completes the proof.\qed

\begin{center}
{ \large Acknowledgements }\end{center} The result of this work
forms a  part of the author's Ph.D. thesis at the University of
Isfahan. The author thanks her supervisor, Professor Saeid Azam,
for his continuous guidance and encouragement. The author also
would like to thank Professor Erhard Neher for some helpful
discussion during the time when she was a visiting researcher at
the University of Ottawa. She wishes to thank the referee of the
paper for his/her invaluable suggestions.

\end{document}